\documentclass{amsart}

\usepackage{rlepsf,amssymb}

\newtheorem{thm}{Theorem}[section]
\newtheorem{lem}[thm]{Lemma}
\newtheorem{cor}[thm]{Corollary}
\newtheorem{prop}[thm]{Proposition}
\newtheorem{rem}[thm]{Remark}

\newtheorem*{ttm}{Theorem}
\newtheorem*{fct}{Fact}
\newtheorem{quest}{Question}
\newtheorem{defn}[thm]{Definition}

\newtheorem*{claim}{Claim}
\newtheoremstyle{definition}{7pt plus6.3pt minus6.3pt}{7pt plus3pt minus3pt}%
{\rm}{}{\bf}{}{0.75em}{\thmname{#1}\thmnumber{ #2}\thmnote{\sl\stdspace#3}}
\theoremstyle{definition}\newtheorem{example}[thm]{Example}
\newtheorem{exercise}[thm]{\small Exercise}

\newcommand{\bbr}{\begin{rem}\em} 
\newcommand{\eer}{\end{rem}}

\newcommand{\bex}{\begin{example}} 
\newcommand{\eex}{\end{example}}

\newcommand{\bhw}{\begin{exercise}\small} 
\newcommand{\ehw}{\end{exercise}}

\newcommand{\be}{\begin{enumerate}}
\newcommand{\ee}{\end{enumerate}}

\def\C{\hbox{$\mathbb C$} }

\def\R{\hbox{$\mathbb R$} }

\def\dfn#1{{\em #1}}

\begin{document}

\title{Lectures on open book decompositions and contact structures}

\author{John B. Etnyre}
\address{University of Pennsylvania, Philadelphia, PA 19104}
\email{etnyre@math.upenn.edu}
\urladdr{http://www.math.upenn.edu/\char126 etnyre}

\thanks{From lecturers given at the 2004 Clay Mathematics Institute Summer School on Floer Homology, 
Gauge Theory, and Low Dimensional Topology at the Alfr\'ed R\'enyi Institute; 
www.claymath.org/programs/summer\_school/2004/.}


\maketitle


\section{Introduction}
The main goal of this survey is to discuss the proof and examine some consequences of the following
fundamental theorem of Giroux.
\begin{thm}[Giroux 2000, \cite{Giroux??}]\label{thm:main}
Let $M$ be a closed oriented 3-manifold. Then there is a one to one correspondence between
\[\{\text{oriented contact structures on } M \text{ up to isotopy}\}\]
and
\[\{\text{open book decompositions of } M \text{ up to positive stabilization}\}.\]
\end{thm}
This theorem plays a pivotal role in studying cobordisms of contact structures and understanding
filling properties of contact structures, see \cite{AO, Eliashberg04,  Etnyre??, Etnyre04,
EtnyreHonda02a, Gay02}. This better understanding of fillings leads to various topological
applications of contact geometry. Specifically, the much studied property P for knots was
established by P.~Kronheimer and T.~Mrowka in \cite{KronheimerMrowka}. A non-trivial knot has property P
if non-trivial surgery on it never gives a homotopy sphere. In addition P.~Ozsv\'ath and Z.~Szab\'o in
\cite{OzsvathSzabo} gave an alternate proof of a characterization of the unknot via surgery which was
originally established in \cite{KMOS}. This characterization says that the unknot is the only knot
on which $p$-surgery yields $-L(p,1).$ Moreover, in \cite{OzsvathSzabo} it is shown that the Thurston
norm is determined by Heegaard Floer Homology.

Ideally the reader should be familiar with low-dimensional topology at the level of, say
\cite{Rolfsen}. In particular, we will assume familiarity with Dehn surgery, mapping tori and basic
algebraic topology. At various points we also discuss branch coverings, Heegaard splittings and
other notions; however, the reader unfamiliar with these notions should be able to skim these parts
of the paper without missing much, if any, of the main line of the arguments. Since diffeomorphisms
of surfaces play a central role in much of the paper and specific conventions are important we have
included an Appendix discussing basic facts about this. We also assume the reader has some
familiarity with contact geometry. Having read \cite{EtnyreCIntro} should be sufficient background
for this paper. In order to accommodate the reader with little background in contact geometry we
have included brief discussions, scattered throughout the paper, of all the necessary facts. Other
good introductions to contact geometry are \cite{Abetal, GeigesIntro}, though a basic understanding
of convex surfaces is also useful but is not covered in these sources.

In the next three sections we give a thorough sketch of the proof of Theorem~\ref{thm:main}. In
Section~\ref{sec:ob} we define open book decompositions of 3-manifolds, discuss their existence and
various constructions. The following two sections discuss how to get a contact structure from an
open book and an open book from a contact structure, respectively. Finally in Section~\ref{sec:apps}
we will consider various applications of Theorem~\ref{thm:main}. While we prove various things about
open books and contact structures our main goal is to prove the following theorem which is the basis
for most of the above mentioned applications of contact geometry to topology.
\begin{ttm}[Eliashberg 2004 \cite{Eliashberg04}; Etnyre 2004 \cite{Etnyre04}]
If $(X,\omega)$ is a symplectic filling of $(M,\xi)$ then there is a closed symplectic manifold
$(W,\omega')$ and a symplectic embedding $(X,\omega)\to (W,\omega').$
\end{ttm}

Acknowledgments: I am grateful to David Alexandre Ellwood, Peter Ozsv\'ath, Andr\'as Stipsicz,
Zoltan Szab\'o, the Clay Mathematics Institute and the Alfr\'ed R\'enyi Institute of Mathematics for
organizing the excellent summer school on ``Floer Homology, Gauge Theory, and Low Dimensional
Topology'' and for giving me an opportunity to give the lectures on which these notes are based. I
also thank Emmanuel Giroux who gave a beautiful series of lectures at Stanford University in 2000
where I was first exposed to the strong relation between open books and contact structures. I am
also grateful to Noah Goodman, Gordana Mati\'c, Andr\'as N\'emethi and Burak Ozbagci for many illuminating
conversations. Finally I thank Paolo Lisca, Stephan Schoenenberger and the referee for valuable comments on the 
first draft of this paper.  This work was supported in part by NSF CAREER Grant (DMS--0239600) and FRG-0244663.

\section{Open book decompositions of 3-manifolds}\label{sec:ob}

Throughout this section (and these notes)
\[\text{\em $M$ is always a closed oriented 3-manifold.}\]
We also mention that when inducing an orientation on the boundary of a manifold we use the ``outward
normal first'' convention. That is, given an oriented manifold $N$ then $v_1,\ldots, v_{n-1}$ is an
oriented basis for $\partial N$ if $\nu, v_1,\ldots, v_{n-1}$ is an oriented basis for $N.$

\begin{defn}
An \dfn{open book decomposition of $M$} is a pair $(B,\pi)$ where
\begin{enumerate}
\item $B$ is an oriented link in $M$ called the \dfn{binding} of the open book and
\item $\pi: M\setminus B \to S^1$ is a fibration of the complement of $B$ such that
  $\pi^{-1}(\theta)$ is the interior of a compact surface $\Sigma_\theta\subset M$ and $\partial
  \Sigma_\theta=B$ for all $\theta\in S^1.$ The surface $\Sigma=\Sigma_\theta,$ for any $\theta,$ is
  called the \dfn{page} of the open book.
\end{enumerate}
\end{defn}
One should note that it is important to include the projection in the data for an open book, since
$B$ does not determine the open book, as the following example shows.
\bex
Let $M=S^1\times S^2$ and $B=S^1\times\{N,S\},$ where $N,S\in S^2.$ There are many ways to fiber
$M\setminus B=S^1\times S^1\times[0,1].$ In particular if $\gamma_n$ is an embedded curve on $T^1$
in the homology class $(1,n),$ then $M\setminus B$ can be fibered by annuli parallel to
$\gamma_n\times[0,1].$ There are diffeomorphisms of $S^1\times S^2$ that relate all of these
fibrations but the fibrations coming from $\gamma_0$ and $\gamma_1$ are not isotopic. There are
examples of fibrations that are not even diffeomorphic.
\eex
\begin{defn}
An \dfn{abstract open book} is a pair $(\Sigma,\phi)$ where
\begin{enumerate}
\item $\Sigma$ is an oriented compact surface with boundary and
\item $\phi:\Sigma\to \Sigma$ is a diffeomorphism such that $\phi$ is the identity
  in a neighborhood of $\partial \Sigma.$ The map $\phi$ is called the \dfn{monodromy}.
\end{enumerate}
\end{defn}

We begin by observing that given an abstract open book $(\Sigma,\phi)$ we get a 3-manifold $M_\phi$
as follows:
\[M_\phi=\Sigma_\phi \cup_\psi \left(\coprod_{|\partial \Sigma|} S^1\times D^2\right),\]
where $|\partial \Sigma|$ denotes the number of boundary components of $\Sigma$ and $\Sigma_\phi$
is the mapping torus of $\phi.$ By this we mean
\[
\Sigma\times[0,1]/\sim,
\]
where $\sim$ is the equivalence relation $(\phi(x),0)\sim(x,1)$ for all $x\in\Sigma.$ Finally,
$\cup_\psi$ means that the diffeomorphism $\psi$ is used to identify the boundaries of the 
two manifolds. For each boundary component $l$ of $\Sigma$ the map $\psi:\partial
(S^1\times D^2) \to l\times S^1\subset \Sigma_\phi$ is defined to be the unique (up to isotopy)
diffeomorphism that takes $S^1\times \{p\}$ to $l$ where $p\in \partial D^2$ and $\{q\}\times
\partial D^2$ to $(\{q'\}\times [0,1]/\sim)=S^1,$ where $q\in S^1$ and $q'\in\partial \Sigma.$ We
denote the cores of the solid tori in the definition of $M_\phi$ by $B_\phi.$

Two abstract open book decompositions $(\Sigma_1,\phi_1)$ and $(\Sigma_2,\phi_2)$ are called
\dfn{equivalent} if there is a diffeomrophism $h:\Sigma_1\to \Sigma_2$ such that $h\circ \phi_2 =
\phi_1\circ h.$

\begin{lem}
We have the following basic facts about open books and abstract open books:
\begin{enumerate}
\item An open book decomposition $(B,\pi)$ of $M$ gives an abstract open book
  $(\Sigma_\pi,\phi_\pi)$ such that $(M_{\phi_\pi}, B_{\phi_\pi})$ is diffeomorphic to $(M,B).$
\item An abstract open book determines $M_\phi$ and an open book $(B_\phi, \pi_\phi)$ up to
diffeomorphism.
\item Equivalent open books give diffeomorphic 3-manifolds.
\end{enumerate}
\end{lem}
\bhw
Prove this lemma.
\ehw
\bbr
Clearly the two notions of open book decomposition are closely related. The basic difference is that when
discussing open books (non-abstract) we can discuss the binding and pages up to {\em isotopy}
in $M,$ whereas when discussing abstract open books we can only discuss them up to {\em
diffeomorphism}. Thus when discussing Giroux's Theorem~\ref{thm:main} we need to use (non-abstract)
open books; however, it is still quite useful to consider abstract open books and we will frequently
not make much of a distinction between them.
\eer
\bex
Let $S^3$ be the unit sphere in $\C^2,$ and $(z_1,z_2)=(r_1e^{i\theta_1}, r_2e^{i\theta_2})$ be
coordinates on $\C^2.$
\begin{enumerate}
\item Let $U=\{z_1=0\}=\{r_1=0\}\subset S^3.$ Thus $U$ is a unit $S^1$ sitting in $S^3.$ It is easy
to see that $U$ is an unknotted $S^1$ in $S^3.$ The complement of $U$ fibers:
\[\pi_U:S^3\setminus U\to S^1: (z_1,z_2)\mapsto \frac{z_1}{|z_1|}.\]
In polar coordinates this map is just $\pi_U(r_1e^{i\theta_1}, r_2e^{i\theta_2})=\theta_1.$ This
fibration is related to the well known fact that $S^3$ is the union of two solid tori. Pictorially
we see this fibration in Figure~\ref{fig:S3}.
\begin{figure}[ht]
  \relabelbox \small {\epsfysize=1.6in\centerline{\epsfbox{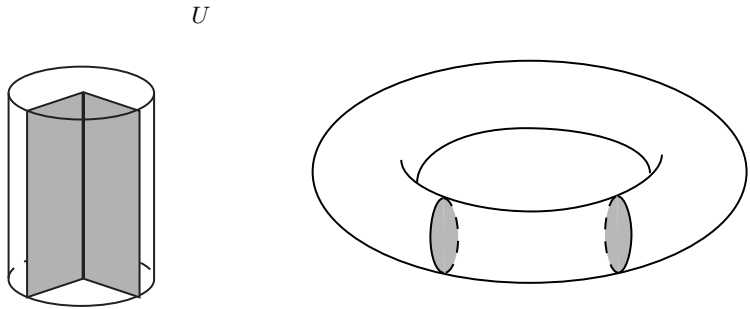}}} \relabel{U}{$U$}
  \endrelabelbox
        \caption{$S^3$ broken into two solid tori (to get the one on the left identify top and bottom 
          of the cylinder). The union of the shaded annuli and disks give two pages in the open book.}
        \label{fig:S3}
\end{figure}
\item Let $H^+=\{(z_1,z_2)\in S^3 : z_1z_2=0\}$ and $H^-=\{(z_1,z_2)\in S^3: z_1\overline{z_2}=0\}.$
\bhw
Show $H^+$ is the positive Hopf link and $H^-$ is the negative Hopf link. See Figure~\ref{fig:hopf}.
(Recall $H^+$ gets an orientation as the boundary of a complex hypersurface in $\C^2,$ and $H^-$ may
be similarly oriented.)
\begin{figure}[ht]
  \relabelbox \small {\epsfysize=1.5in\centerline{\epsfbox{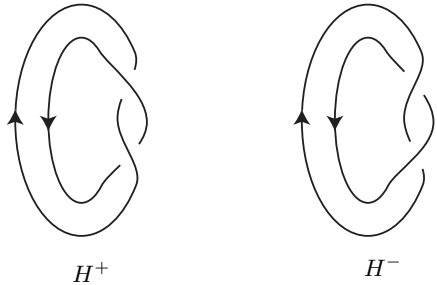}}} 
  \relabel{P}{$H^+$}
  \relabel{M}{$H^-$}
  \endrelabelbox
        \caption{The two Hopf links.}
        \label{fig:hopf}
\end{figure}
\ehw
We have the fibrations
\[
\pi_+:S^3\setminus H^+\to S^1: (z_1,z_2)\mapsto \frac{z_1z_2}{|z_1z_2|}, \text{and}
\]
\[
\pi_-:S^3\setminus H^-\to S^1: (z_1,z_2)\mapsto \frac{z_1\overline{z_2}}{|z_1\overline{z_2}|}.
\]
In polar coordinates these maps are just $\pi_\pm(r_1e^{i\theta_1}, r_2e^{i\theta_2})=\theta_1\pm 
\theta_2.$
\bhw
Picture these fibrations.
\ehw

\item More generally, let $f:\C^2\to \C$ be a polynomial that vanishes at $(0,0)$ and has no
critical points inside $S^3$ except possibly $(0,0).$ Then $B=f^{-1}(0)\cap S^3$ gives an open book
of $S^3$ with fibration
\[
\pi_f: S^3\setminus B \to S^1: (z_1,z_2) \mapsto \frac{f(z_1,z_2)}{|f(z_1,z_2)|}.
\]
This is called the Milnor fibration of the hypersurface singularity $(0,0)$ (note that $(0,0)$ does not
have to be a singularity, but if it is not then $B$ is always the unknot). See \cite{Milnor}.
\end{enumerate}
\eex
\bhw\label{trivmongivesconsum}
Suppose $\Sigma$ is a surface of genus $g$ with $n$ boundary components and $\phi$ is the identity 
map on $\Sigma.$ Show $M_\phi=\#_{2g+n-1} S^1\times S^2.$\hfill\break
HINT: If $a$ is a properly embedded arc in $\Sigma$ then $a\times[0,1]$ is an annulus in the mapping torus $\Sigma_\phi$ that
can be capped off into a sphere using two disks in the neighborhood of the binding.
\ehw

\begin{figure}[ht]
  \relabelbox \small {\epsfysize=1.2in\centerline{\epsfbox{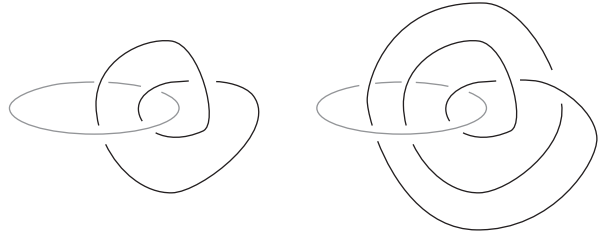}}} 
  \endrelabelbox
        \caption{Here are two links each of which is a link of unknots. If we do a cyclic branched cover of 
          $S^3$ over the grey component and lift the black component to the cover it will become 
          the binding of an open book decomposition of $S^3.$}
        \label{fig:cover}
\end{figure}

\begin{thm}[Alexander 1920, \cite{Alexander20}]\label{openbookexist}
Every closed oriented 3-manifold has an open book decomposition.
\end{thm}
We will sketch three proofs of this theorem.
\begin{proof}[First Sketch of Proof]
We first need two facts
\begin{fct}[Alexander 1920, \cite{Alexander20}]
Every closed oriented 3-manifold $M$ is a branched cover of $S^3$ with branched set some link $L_M.$
\end{fct}
\begin{fct}[Alexander 1923, \cite{Alexander23}]
Every link $L$ in $S^3$ can be braided about the unknot.
\end{fct}
When we say $L$ can be \dfn{braided about the unknot} we mean that if $S^1\times D^2=S^3\setminus U$
then we can isotop $L$ so that $L\subset S^1\times D^2$ and $L$ is transverse to $\{p\}\times D^2$
for all $p\in S^1.$

Now given $M$ and $L_M\subset S^3$ as in the first fact we can braid $L_M$ about the unknot $U.$ Let
$P:M\to S^3$ be the branch covering map. Set $B=P^{-1}(U)\subset M.$ We claim that $B$ is the
binding of an open book. The fibering of the complement of $B$ is simply $\pi=\pi_U\circ P,$ where
$\pi_U$ is the fibering of the complement of $U$ in $S^3.$
\bhw
Prove this last assertion and try to picture the fibration.
\ehw
\end{proof}
Before we continue with our two other proofs let's have some fun with branched covers.
\bhw
Use the branched covering idea in the previous proof to find various open books of
$S^3.$\hfill\break 
HINT: Any cyclic branched cover of $S^3$ over the unknot is $S^3.$ Consider Figure~\ref{fig:cover}.
See also \cite{Goldsmith, Rolfsen}.
\ehw

\begin{proof}[Second Sketch of Proof]
This proof comes from Rolfsen's book \cite{Rolfsen} and relies on the following fact.
\begin{fct}[Lickorish 1962, \cite{Lickorish62}; Wallace 1960, \cite{Wallace60}]
Every closed oriented 3-manifold may be obtained by $\pm 1$ surgery on a link $L_M$ of unknots.
Moreover, there is an unknot $U$ such that $L_M$ is braided about $U$ and each component of $L_M$
can be assumed to link $U$ trivially one time. See Figure~\ref{fig:temp}.
\begin{figure}[ht]
  \relabelbox \small {\epsfysize=1.5in\centerline{\epsfbox{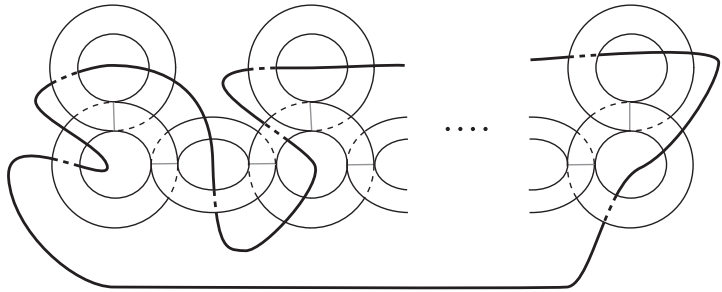}}} \endrelabelbox
        \caption{All the unknots in the link $L_M$ can be isotoped to be on the annuli depicted 
          here. The heavy black line is the unknot $U.$}
        \label{fig:temp}
\end{figure}
\end{fct}
Now $(U,\pi_U)$ is an open book for $S^3.$ Let $N$ be a small tubular neighborhood of $L_M.$ Each
component $N_C$ of $N$ corresponds to a component $C$ of $L_M$ and we can assume that $N_C$
intersects the fibers of the fibration $\pi_U$ in meridional disks. So the complement of $U\cup N$
fibers so that each $\partial N_C$ is fibered by meridional circles. To perform $\pm 1$ surgery on
$L_M$ we remove each of the $N_C$'s and glue it back sending the boundary of the meridional disk to
a $(1,\pm 1)$ curve on the appropriate boundary component of $\overline {S^3\setminus N}.$ After the
surgery we have $M$ and inside $M$ we have the union of surgery tori $N',$ the components of which
we denote $N_C'$ and the cores of which we denote $C'.$ We denote the union of the cores by $L'.$
Inside $M$ we also have the ``unknot'' $U$ (of course $U$ may not be an unknot any more, for example
it could represent non-trivial homology in $M$). Since $M\setminus (U\cup N') =S^3\setminus (U\cup
N)$ we have a fibration of $M\setminus (U\cup N')$ and it is easy to see that the fibration induces on
$\partial N'_C$ a fibration by $(1,\pm 1)$ curves. We can fiber $N'_C\setminus C'$ by annuli so that
the induced fibration on $\partial N'_C$ is by $(1,\pm 1)$ curves. Thus we may extend the fibration
of $M\setminus (U\cup N')$ to a fibration of $M\setminus(U\cup L'),$ hence inducing an open book of
$M.$
\bhw
Convince yourself of these last statements.
\ehw
\bbr
We have produced an open book for $M$ with planar pages!
\eer
\end{proof}

\begin{proof}[Third Sketch of Proof]
This proof is due to Harer. We need 
\begin{fct}[Harer 1979, \cite{Harer}]
An oriented compact 4-manifold has an achiral Lefschetz fibration with non-closed leaves over a disk if and
only if it admits a handle decomposition with only 0-, 1-, and 2-handles.
\end{fct}
An achiral Lefschetz fibration of a 4-manifold $X$ over a surface $S$ is simply a map $\pi:X\to S$
such that the differential $d\pi$ is onto for all but a finite number of points $p_1,\ldots p_k \in
\text{int}(X),$ where there are complex coordinate charts $U_i$ of $p_i$ and $V_i$ of $\pi(p_i)$ such
that $\pi_{U_i}(z_1, z_2) = z_1^2+z_2^2.$ Note the definition implies that $\pi$ restricted to
$X\setminus \pi^{-1}(\pi(\{p_1\ldots, p_k\}))$ is a locally trivial fibration. We denote a generic
fiber by $\Sigma_\pi.$
\begin{fct}[Lickorish 1962, \cite{Lickorish62}; Wallace 1960, \cite{Wallace60}]
Every closed oriented 3-manifold is the boundary of a 4-manifold built with only 0- and 2-handles.
\end{fct}
Given a 3-manifold $M$ we use this fact to find a 4-manifold $X$ with $\partial X=M$ and $X$ built
with only 0- and 2-handles. Then the previous fact gives us an achiral Lefschetz fibration $\pi:X\to
D^2.$ Set $B=\partial \pi^{-1}(x)$ for a non-critical value $x\in \text{int}(D^2).$ We claim that
$B$ is the binding of an open book decomposition for $M$ and the fibration of the complement is the
restriction of $\pi$ to $M\setminus B.$
\end{proof}

\begin{defn}
Given two abstract open books $(\Sigma_i,\phi_i), i=0,1,$ let $c_i$ be an arc properly embedded in
$\Sigma_i$ and $R_i$ a rectangular neighborhood of $c_i,$ $R_i=c_i\times[-1,1].$ The \dfn{Murasugi
sum} of $(\Sigma_0,\phi_0)$ and $(\Sigma_1,\phi_1)$ is the open book
$(\Sigma_0,\phi_0)*(\Sigma_1,\phi_1)$ with page
\[\Sigma_0*\Sigma_1= \Sigma_0 \cup _{R_1=R_2} \Sigma_1,\]
where $R_0$ and $R_1$ are identified so that $c_i\times\{-1,1\}=(\partial c_{i+1})\times [-1,1],$
and the monodromy is $\phi_0\circ \phi_1.$
\end{defn}

\begin{thm}[Gabai 1983, \cite{Gabai83}]\label{thm:murasugi}
\[
M_{(\Sigma_0,\phi_0)}\#M_{(\Sigma_1,\phi_1)} \text{ is diffeomorphic to } 
M_{(\Sigma_0,\phi_0)*(\Sigma_1,\phi_1)}.
\]
\end{thm}
\begin{proof}[Sketch of Proof]
The proof is essentially contained in Figure~\ref{fig:msum}. The idea is that
$B_0=R_0\times[\frac12,1]$ is a 3-ball in $M_{(\Sigma_1,\phi_1)}$ and similarly for
$B_1=R_1\times[0,\frac12]$ in $M_{(\Sigma_0,\phi_0)}.$ Now $(\Sigma_0*\Sigma_1)\times[0,1]$ can be
formed as shown in Figure~\ref{fig:msum}.

\begin{figure}[ht]
  \relabelbox \small {\epsfysize=2.5in\centerline{\epsfbox{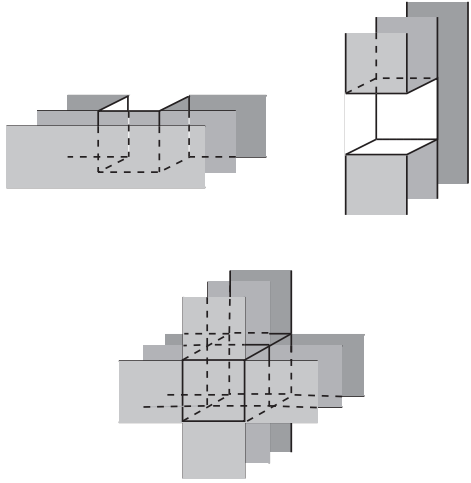}}} \endrelabelbox
        \caption{At the top left is a piece of $\Sigma_0\times[0,1]$ near $c_0$ with $B_0$ cut out. 
          The lightest shaded part is $\Sigma_0\times\{0\}$ the medium shaded part is $\Sigma_0\times 
          \{\frac12\}$ and the darkest shaded part is $\Sigma_0\times\{1\}.$ The top right is a similar 
          picture for $\Sigma_1.$ The bottom picture is $(\Sigma_0*\Sigma_1)\times[0,1].$}
        \label{fig:msum}
\end{figure}
Think about forming the mapping cylinder of $\phi_0$ by gluing $\Sigma_0\times\{0\}$ to
$\Sigma_0\times\{1\}$ using the identity and then cutting the resulting $\Sigma_0\times S^1$ along
$\Sigma_0\times\{\frac14\}$ and regluing using $\phi_0.$ Similarly think about the mapping cylinder
for $\phi_1$ as $\Sigma_1\times S^1$ reglued along $\Sigma_1\times \{\frac34\}$ and the mapping
cylinder for $\phi_0\circ \phi_1$ as $(\Sigma_0*\Sigma_1)\times S^1$ reglued by $\phi_0$ along
$(\Sigma_0*\Sigma_1)\times \{\frac14\}$ and by $\phi_1$ along $(\Sigma_0*\Sigma_1)\times
\{\frac34\}.$ Thus we see how to fit all the mapping cylinders together nicely.
\bhw
Think about how the binding fits in and complete the proof.
\ehw
\end{proof}

\begin{figure}[ht]
  \relabelbox \small {\epsfysize=1.2in\centerline{\epsfbox{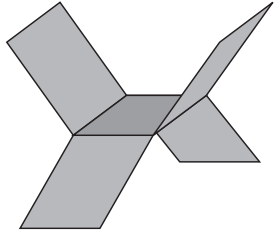}}} 
  \endrelabelbox
        \caption{Two ambient surfaces being summed together.}
        \label{fig:amb}
\end{figure}

\begin{defn}
A \dfn{positive (negative) stabilization} of an abstract open book $(\Sigma,\phi)$ is the open book
\begin{enumerate}
\item with page $\Sigma'=\Sigma\cup \text{ 1-handle}$ and
\item monodromy $\phi'=\phi\circ \tau_c$ where $\tau_c$ is a right- (left-)handed Dehn twist along a
curve $c$ in $\Sigma'$ that intersects the co-core of the 1-handle exactly one time.
\end{enumerate}
We denote this stabilization by $S_{(a,\pm)}(\Sigma,\phi)$ where $a=c\cap \Sigma$ and $\pm$ refers
to the positivity or negativity of the stabilization. (We omit the $a$ if it is unimportant in a
given context.)
\end{defn}
\bhw
Show \[S_\pm(\Sigma,\phi)=(\Sigma,\phi)*(H^\pm,\pi_\pm)\] where $H^\pm$ is the positive/negative 
Hopf link and $\pi_\pm$ is the corresponding fibration of its complement.  
\ehw
From this exercise and Theorem~\ref{thm:murasugi} we immediately have:
\begin{cor}
\[M_{(S_\pm(\Sigma,\phi))}=M_{(\Sigma,\phi)}.\]
\end{cor}
\bhw
Show how to do a Murasugi sum ambiently. That is show how to perform a Murasugi sum for open book 
decompositions (not abstract open books!). Of course one of the open books must be an open book for 
$S^3.$
\hfill\break
HINT: See Figure~\ref{fig:amb}.
\ehw

\bhw
Use Murasugi sums to show the right- and left-handed trefoil knots and the figure eight knot all 
give open book decompositions for $S^3.$
\ehw
\bhw
Use Murasugi sums to show all torus links give open book decompositions of $S^3.$
\ehw
\bhw
Show that every 3-manifold has an open book decomposition with connected binding.
\ehw
\bhw
Use the previous exercise to prove a theorem of Bing: A closed oriented 3-manifold is $S^3$ if 
and only if every simple closed curve 
in $M$ is contained in a 3-ball.
\hfill\break
HINT: The only surface bundle over $S^1$ yielding an orientable manifold that is not irreducible 
is $S^1\times S^2.$
\ehw

\section{From open books to contact structures}
\begin{defn}
An \dfn{(oriented) contact structure} $\xi$ on $M$ is an oriented plane field $\xi\subset TM$ for
which there is a 1-form $\alpha$ such that $\xi=\ker \alpha$ and $\alpha\wedge d\alpha >0.$ (Recall:
$M$ is oriented.)
\end{defn}
\bbr\label{rem:poscont}
What we have really defined is a {\em positive} contact structure, but since this is all we will 
talk about we will stick to this definition. 
\eer
\bex
\begin{enumerate}
\item On $\R^3$ we have the standard contact structure $\xi_{std}=\ker(dz+r^2d\theta).$ See
Figure~\ref{fig:fex}.
\begin{figure}[ht]
  \relabelbox \small {\epsfysize=2.4in\centerline{\epsfbox{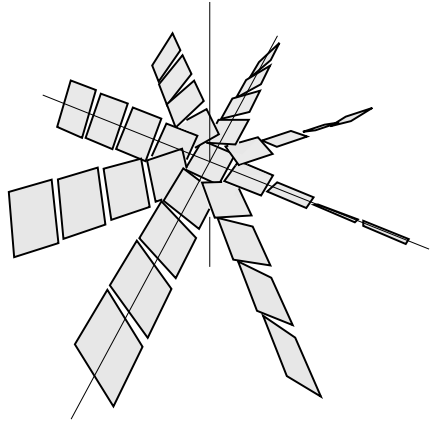}}} \endrelabelbox
        \caption{Standard contact structure on $\R^3.$ (Picture by Stephan Sch\"onenberger.)}
        \label{fig:fex}
\end{figure}
\item On $S^3,$ thought of as the unit sphere in $\C^2,$ we have $\xi_{std}$ the set of complex
tangents. That is $\xi_{std}=TS^3\cap i(TS^3).$ We can also describe this plane field as
$\xi_{std}=\ker(r_1^2d\theta_1+r_2^2d\theta_2),$ where we are using coordinates
$(r_1e^{i\theta_1},r_2e^{i\theta_2})$ on $\C^2.$
\end{enumerate}
\eex

We will need the following facts about contact structures. Most of these facts are proven or
discussed in \cite{Abetal, EtnyreCIntro, GeigesIntro}.
\begin{itemize}
\item All 3-manifolds admit a contact structure. See Theorem~\ref{twthm} below.
\item Locally all contact structures look the same. This is called Darboux's Theorem and means that
if $p_i$ is a point in the contact manifold $(M_i,\xi_i), i=0,1,$ then there is a neighborhood $U_i$
of $p_i$ and a diffeomorphism $f:U_0\to U_1$ such that $f_*(\xi_0)=\xi_1.$ Such a diffeomorphism is
called a \dfn{contactomorphism}.
\item Given two contact manifolds $(M_i,\xi_i), i=0,1,$ we can form their \dfn{contact connected
sum} $(M_0\# M_1, \xi_0 \# \xi_1)$ as follows: there are balls $B_i\subset M_i$ and an orientation
reversing diffeomorphism $f:\partial (\overline{M_0\setminus B_0})\to \partial(\overline{M_1\setminus B_1})$
such that $M_0\# M_1$ is formed by gluing $M_0\setminus B_0$ and $M_1\setminus B_1$ together using $f$ and
\[\xi_0|_{M_0\setminus B_0}\cup \xi_1|_{M_1\setminus B_1}\]
extends to a well defined contact structure on $M_0\# M_1.$ See \cite{Colin97}.
\item Given a 1-parameter family of contact structures $\xi_t, t\in[0,1],$ there is a 1-parameter
family of diffeomorphisms $\phi_t:M\to M$ such that $(\phi_t)_*(\xi_0)=\xi_t.$ This is called Gray's
Theorem.
\end{itemize}
Two contact structures $\xi_0$ and $\xi_1$ are called isotopic if there is a 1-parameter family of
contact structures connecting them.

\begin{defn}
A contact structure $\xi$ on $M$ is \dfn{supported} by an open book decomposition $(B,\pi)$ of $M$
if $\xi$ can be isotoped through contact structures so that there is a contact 1-form $\alpha$ for
$\xi$ such that
\begin{enumerate}
\item $d\alpha$ is a positive area form on each page $\Sigma_\theta$ of the open book and
\item $\alpha>0$ on $B$ (Recall: $B$ and the pages are oriented.)
\end{enumerate}
\end{defn}
\begin{lem}\label{equivsupp}
The following statements are equivalent
\begin{enumerate}
\item[(1)] The contact manifold $(M, \xi)$ is supported by the open book $(B,\pi).$
\item[(2)] $(B,\pi)$ is an open book for $M$ and $\xi$ can be isotoped
to be arbitrarily close (as oriented plane fields), on compact subsets of the pages, to the tangent 
planes to the pages of the open book in such a way that after some point in the isotopy the contact
planes are transverse to $B$ and transverse to the pages of the open book in a fixed neighborhood of
$B.$
\item [(3)] $(B,\pi)$ is an open book for $M$ and there is a Reeb vector field $X$ for a contact
structure isotopic to $\xi$ such that $X$ is (positively) tangent to $B$ and (positively) transverse
to the pages of $\pi.$
\end{enumerate}
\end{lem}
The condition in part (2) of the lemma involving transversality to the pages is to prevent excess twisting near the
binding and may be dispensed with for tight contact structures.

Recall that a vector field $X$ is a Reeb vector field for $\xi$ if it is transverse to $\xi$ and its flow
preserves $\xi.$ This is equivalent to saying there is a contact form $\alpha$ for $\xi$ such that
$\alpha(X)=1$ and $\iota_X d\alpha=0.$ The equivalence of (1) and (2) is supposed to give some
intuition about what it means for a contact structure to be supported by an open book. We do not
actually use (2) anywhere in this paper, but it is interesting to know that ``supported'' can be
defined this way. Similarly condition (3) should be illuminating if you have studied Reeb vector fields
in the past.
\begin{proof}
We begin with the equivalence of (1) and (2). Suppose $(M,\xi)$ is supported by $(B,\pi).$ So
$M\setminus B$ fibers over $S^1.$ Let $d\theta$ be the coordinate on $S^1.$ We also use $d\theta$ to
denote the pullback of $d\theta$ to $M\setminus B,$ that is, for $\pi^*d\theta.$ Near each component
for the binding we can choose coordinates $(\psi, (r,\theta))$ on $N=S^1\times D^2$ in such a way
that $d\theta$ in these coordinates and $\pi^*d\theta$ agree. Choosing the neighborhood $N$ small
enough we can assume that $\alpha(\frac{\partial}{\partial \psi})>0$ (since $\alpha$ is positive on
$B$). Choose an increasing non-negative function $f:[0,\epsilon]\to \R$ that equals $r^2$ near 0 and
1 near $\epsilon,$ where $\epsilon$ is chosen so that $\{(\psi, (r, \theta))| r<\epsilon\}\subset
N.$ Now consider the 1-form $\alpha_R=\alpha +R f(r)d\theta,$ where $R$ is any large constant. (Here
we of course mean that outside the region 
$\{(\psi, (r, \theta))| r<\epsilon\}$ we just take $f$ to be 1.) Note that
$\alpha_R$ is a contact 1-form for all $R>0.$ Indeed
\[
\alpha_R\wedge d\alpha_R = \alpha\wedge d\alpha + Rf d\theta \wedge d\alpha 
+ Rf'\alpha \wedge dr\wedge d\theta.
\]
The first term on the right is clearly positive since $\alpha$ is a contact form. The second term is
also positive since $d\alpha$ is a volume form for the pages, $d\theta$ vanishes on the pages and
is positive on the oriented normals to the pages. Finally the last term is non-negative since $dr\wedge
d\theta$ vanishes on $\frac{\partial}{\partial \psi}$ while $\alpha(\frac{\partial}{\partial
\psi})>0.$ As $R\to \infty$ we have a 1-parameter family of contact structures $\xi_R=\ker
\alpha_R$ that starts at $\xi=\xi_0$ and converges to the pages of the open book away from the
binding while staying transverse to the binding (and the pages near the binding).

Now for the converse we assume (2). Let $\xi_s$ be a family of plane fields isotopic to $\xi$ that
converge to a singular plane field tangent to the pages of the open book (and singular along the
binding) as $s\to \infty.$ Let $\alpha_s$ be contact forms for the $\xi_s.$ We clearly have that
$\alpha_s>0$ on $B.$ 

Thinking of $M\setminus B$ as the mapping torus $\Sigma_\phi$ we can use 
coordinates $(x,\theta)\in \Sigma\times[0,1]$ (we use $\theta$ for the coordinate on $[0,1]$
since on the mapping torus $\Sigma_\phi$ this is the pullback of $\theta$ on $S^1$ under the fibration) and write
\[
\alpha_s = \beta_s(\theta) +u_s(\theta)\, d\theta,
\]
where $\beta_s(\theta)$ is a 1-form on $\Sigma$ and $u_s(\theta)$ is a function on $\Sigma$ for each
$s$ and $\theta.$ Let $N$ be a tubular neighborhood of $B$ on which $\xi_s$ is transverse to the pages of
the open book and let $N'$ be a tubular neighborhood of $B$ contained in $N.$ We can choose $N'$ so that 
$\overline{M\setminus N'}$ is a mapping torus $\Sigma'_{\phi'}$ where $\Sigma'\subset \Sigma$ is $\Sigma$ minus
a collar neighborhood of $\partial \Sigma$ and $\phi'=\phi|_{\Sigma'}$ is the identity near $\partial \Sigma'.$
For $s$ large enough $u_s(\theta)>0$ on $\overline{M\setminus N'}$ for all $\theta,$ since $\alpha_s$ converges
uniformly to some positive multiple of $d\theta$ on $\overline{M\setminus N'}$ as $s\to \infty.$
Thus, for large $s,$ on $\overline{M\setminus N'}=\Sigma'_{\phi'}$ we can divide
$\alpha_s$ by $u_s(\theta)$ and get a new family of contact forms
\[
\alpha'_s= \beta'_s(\theta) + d\theta.
\]
We now claim
that $d\alpha'_s|_{\text{page}}=d\beta'_s$ is a positive volume form on $\Sigma'.$ To see this
note that
\[
\alpha'_s\wedge d\alpha'_s = d\theta \wedge (d\beta_s'(\theta)-\beta'_s(\theta)\wedge \frac{\partial
\beta'_s}{\partial \theta}(\theta)).
\]
So clearly
\begin{equation}\label{eq:equiv1}
d\beta_s'(\theta)-\beta'_s(\theta)\wedge \frac{\partial \beta'_s}{\partial \theta}(\theta)>0.
\end{equation} 
To see that $d\beta'_s(\theta)>0$ for $s$ large enough, we note that the second term in this equation
vanishes to higher order than the first as $s$ goes to infinity. From this one can easily conclude
that $d\beta'_s(\theta)\geq 0$ for $s$ large enough. 
\bhw
Verify this last statement.\hfill\break
HINT: Assume, with out loss of generality, the 1-forms $\alpha_s$ are analytic and are analytic in $s.$ 
\ehw
By adding a small multiple of a 1-form, similar to the one constructed on the mapping torus in the proof of 
Theorem~\ref{twthm} below,
we easily see that for a fixed $s,$ large enough, we can assume $d\beta'_s(\theta)>0.$
\bhw
Show how to add this 1-form to $\alpha'_s$ preserving all the properties of $\alpha'_s$ in $N'$
but still having $d\alpha'_s>0$ on $\Sigma'_\theta.$\hfill\break
HINT: Make sure the 1-form and its derivative are very small and use a cutoff function that is $C^1$-small too. 
\ehw 

We may now assume that 
$d\alpha'_s$ is a volume form on the pages of $\Sigma'_{\phi'}.$ Denote $\alpha'_s$ by $\alpha.$
We are left to verify $\alpha$ can be modified to have the desired properties in $N.$ 
\bhw
Show that we may assume each component of $N$ is diffeomorphic to $S^1\times D^2$ with
coordinates $(\psi, (r, \theta))$ such that the pages of the open book go to constant $\theta$
annuli in $S^1\times D^2$ and the contact structure $\ker \alpha$ on $M$ maps to $\ker (d\psi+f(r)\, d\theta),$
for some function $f(r).$
\hfill\break 
HINT: This is more than a standard neighborhood of a transverse curve. Think about the foliation on 
the pages of the open book near the binding and on the constant $\theta$ annuli. 
\ehw
Under this identification $\alpha$ maps to some contact form $\alpha'=h\,  (d\psi+f(r)\, d\theta)$
near the boundary of $S^1\times D^2,$ where $h$ is function on this neighborhood. By scaling
$\alpha$ if necessary we may assume that $h>1$ where it is defined.
\bhw
Show that $d\alpha'$ is a volume form on the (parts of the) constant $\theta$ annuli (where it is 
defined) if and only if $h_r>0.$
\ehw
Since that we know $d\alpha$ is a volume form on the pages of the mapping torus, $h_r>0$ where it is
defined. Moreover we can
extend it to all of $S^1\times D^2$ so that it is equal to 1 on $r=0$ and so that $h_r>0$
everywhere. Thus the contact form equal to $\alpha$ off of $N'$ and
equal to $h(d\psi+f(r)\, d\theta)$ on each component of $N$ is a globally defined contact form for
$\xi_s$ and satisfies conditions (1) of the lemma.
\bhw
Try to show that the condition in part (2) of the lemma involving transversality to the pages
of the open book near $B$ is unnecessary if the contact structure is tight.\hfill\break
HINT: There is a unique universally tight contact structure on a solid torus with a fixed non-singular
characteristic foliation on the boundary that is transverse to the meridional circles.
\ehw

We now establish the equivalence of (1) and (3). Assume (3) and let $X$ be the vector field
discussed in (3). Since $X$ is positively tangent to the binding we have $\alpha>0$ on oriented
tangent vectors to $B.$ Moreover, since $X$ is positively transverse to the pages of the open book
we have $d\alpha=\iota_X \alpha\wedge d\alpha>0$ on the pages. Thus $(M,\xi)$ is supported by
$(B,\pi).$ Conversely assume (1) is true and let $\alpha$ be the contact form implicated in the
definition of supporting open book. Let $X$ be the Reeb vector field associated to $\alpha.$ It is
clear that $X$ is positively transverse to the pages of the open book since $d\alpha$ is a volume
form on the pages. Thus we are left to check that $X$ is positively tangent to $B.$ To this end
consider coordinates $(\psi, (r,\theta))$ on a neighborhood of a component of $B$ such that constant
$\theta$'s give the pages of the open book in the neighborhood. Switching $(r,\theta)$ coordinates
to Cartesian coordinates $(x,y)$ we can write $X=f\frac{\partial}{\partial \psi}
+g\frac{\partial}{\partial x}+h\frac{\partial}{\partial y},$ where $f,g,h$ are functions. We need to
see that $g$ and $h$ are zero when $(x,y)=(0,0).$ This is clear, for if say $g>0$ at some point
$(c,(0,0))$ then it will be positive in some neighborhood of this point in particular at $(c,(0,\pm
\epsilon))$ for sufficiently small $\epsilon.$ But at $(c,(0,\epsilon))$ the
$\frac{\partial}{\partial x}$ component of $X$ must be negative, not positive, in order to be
positively transverse to the pages. Thus $g$ and $h$ are indeed zero along the binding.
\end{proof}
\bex
Let $(U,\pi_U)$ be the open book for $S^3,$ where $U$ is the unknot and 
\[
\pi_U:S^3\setminus U \to S^1: (r_1,\theta_1, r_2,\theta_2) \mapsto \theta_1.
\]
(Recall that we are thinking of $S^3$ as the unit sphere in $\C^2$.) This open book supports the standard
contact structure $\xi_{std}=\ker(r_1^2d \theta_1+r_2^2d\theta_2).$ To see this notice that for fixed
$\theta_1$ the page $\pi_U^{-1}(\theta_1)$ is parameterized by
\[
f(r,\theta)=(\sqrt{1-r^2},\theta_1, r, \theta).
\]
Thus $df^*(r_1^2d\theta_1+r_2^2d\theta_2)=2r\, dr\wedge d\theta$ which is the volume form on the disk.
Moreover the positively oriented tangent to $U$ is $\frac{\partial}{\partial \theta_2}$ and
$\alpha(\frac{\partial}{\partial \theta_2})>0$
\eex
\bhw
Show that $(H^+,\pi_+)$ also supports $\xi_{std}$ but that $(H^-,\pi_-)$ does not.
\ehw

\begin{thm}[Thurston-Winkelnkemper 1975, \cite{TW}]\label{twthm}
Every open book decomposition $(\Sigma,\phi)$ supports a contact structure $\xi_\phi$ on $M_\phi.$
\end{thm}
\begin{proof}
Recall
\[
M_\phi=\Sigma_\phi \cup_\psi \left(\coprod_{|\partial \Sigma|} S^1\times D^2\right),
\]
where $\Sigma_\phi$ is the mapping torus of $\phi.$ We first construct a contact structure on
$\Sigma_\phi.$ To this end we consider the set
\[
\begin{aligned}S=\{ \text{1-forms } \lambda : &
 \text{ (1) }\lambda= (1+s)d\theta \text{ near } \partial \Sigma \text{ and}\\ &\text{ (2) }d\lambda
\text{ is a volume form on } \Sigma\}\end{aligned}
\]
where near each boundary component of $\Sigma$ we use coordinates $(s,\theta)\in [0,1]\times S^1.$
\bhw Show this set is convex. \ehw To show this set is non-empty let $\lambda_1$ be any 1-form on
$\Sigma$ that has the right form near the boundary. Note that
\[\int_\Sigma d\lambda_1 = \int_{\partial \Sigma} \lambda_1 =  2\pi|\partial \Sigma|.\]
Let $\omega$ be any volume form on $\Sigma$ whose integral over $\Sigma$ is $2\pi|\partial \Sigma|$
and near the boundary of $\Sigma$ equals $ds\wedge d\theta.$ We clearly have
\[\int_\Sigma \big( \omega-d\lambda_1\big) =0\]
and $\omega-d\lambda_1=0$ near the boundary. Thus the de Rham theorem says we can find a 1-form
$\beta$ vanishing near the boundary such that $d\beta=\omega-d\lambda_1.$ One may check
$\lambda=\lambda_1+\beta$ is a form in $S.$

Now given $\lambda\in S$ note that $\phi^*\lambda$ is also in $S.$ Consider the 1-form
\[
\lambda_{(t,x)}= t\lambda_x +(1-t)(\phi^*\lambda)_x
\]
on $\Sigma\times[0,1]$ where $(x,t)\in\Sigma\times[0,1]$ and set
\[\alpha_K=\lambda_{(t,x)}+ Kdt.\]
\bhw
Show that for sufficiently large $K$ this form is a contact form.
\ehw
It is clear that this form descends to a contact form on the mapping torus $\Sigma_\phi.$ We now
want to extend this form over the solid tori neighborhood of the binding. To this end consider the
map $\psi$ that glues the solid tori to the mapping torus. In coordinates $(\varphi,(r,\vartheta))$
on $S^1\times D^2$ where $D^2$ is the unit disk in the $\R^2$ with polar coordinates we have
\[
\psi(\varphi, r, \vartheta) = (r-1+\epsilon, -\varphi, \vartheta).
\]
This is a map defined near the boundary of $S^1\times D^2.$ Pulling back the contact form $\alpha_K$
using this maps gives
\[
\alpha_\psi=K\, d\vartheta - (r+\epsilon)\, d\varphi.
\]
We need to extend this over all of $S^1\times D^2.$ We will extend using a form of the form
\[
f(r)\, d\varphi + g(r)\, d\vartheta.
\]
\bhw
Show this form is a contact form if and only if $f(r)g'(r)-f'(r)g(r)>0.$ Said another way, that
\[\begin{pmatrix}
f(r)\\ g(r)
\end{pmatrix},
\begin{pmatrix}
f'(r) \\ g'(r)
\end{pmatrix}
\]
is an oriented basis for $\R^2$ for all $r.$
\ehw
Near the boundary $\alpha_\psi$ is defined with $f(r)= -(r+\epsilon)$ and $g(r)=K.$ Near the core of
$S^1\times D^2$ we would like $f(r)=1$ and $g(r)=r^2.$
\bhw
Show that $f(r)$ and $g(r)$ can be chosen to extend $\alpha_\psi$ across the solid torus. \hfill\break
HINT: Consider the parameterized curve $(f(r), g(r)).$ This curve is defined for $r$ near 0 and 1;
can we extend it over all of $[0,1]$ so that the position and tangent vector are never collinear?
\ehw
\end{proof}

\begin{prop}[Giroux 2000, \cite{Giroux??}]
Two contact structures supported by the same open book are isotopic.
\end{prop}
\begin{proof}
Let $\alpha_0$ and $\alpha_1$ be the contact forms for $\xi_0$ and $\xi_1,$ two contact structures
that are supported by $(B,\pi).$ In the proof of Lemma~\ref{equivsupp} we constructed a contact form
$\alpha_R=\alpha+Rf(r)d\theta$ from $\alpha.$ (See the proof of the lemma for the definitions of the
various terms.) In a similar fashion we can construct $\alpha_{0R}$ and $\alpha_{1R}$ from
$\alpha_0$ and $\alpha_1.$ These are all contact forms for all $R\geq0.$ Now consider
\[
\alpha_s = s \alpha_{1R} + (1-s) \alpha_{0R}.
\]
\bhw
For large $R$ verify that $\alpha_s$ is a contact form for all $0\leq s\leq 1.$ \hfill\break
HINT: There are three regions to consider when verifying that $\alpha_s$ is a contact form. The
region near the binding where $f(t)=r^2,$ the region where $f$ is not 1 and the region where $f$ is
1. Referring back to the proof of Lemma~\ref{equivsupp} should help if you are having difficulty when
considering any of these regions.
\ehw
Thus we have an isotopy from $\alpha_0$ to $\alpha_1.$
\end{proof}
We now know that for each open book $(B,\pi)$ there is a unique contact structure supported by
$(B,\pi).$ We denote this contact structure by $\xi_{(B,\pi)}.$ If we are concerned with the
abstract open books we denote the contact structure by $\xi_{(\Sigma,\phi)}.$
\begin{thm}
We have
\[
\xi_{(\Sigma_0,\phi_0)}\# \xi_{(\Sigma_1,\phi_1)}=\xi_{(\Sigma_0,\phi_0)*(\Sigma_1,\phi_1)}
\]
\qed
\end{thm}
This theorem follows immediately from the proof of Theorem~\ref{thm:murasugi} concerning the effect
of Murasugi sums on the 3--manifold. The theorem seems to have been known is some form or another for
some time now but the first reference in the literature is in Torisu's paper \cite{Torisu00}. 
\bhw
Go back through the proof of Theorem~\ref{thm:murasugi} and verify that the contact structures are 
also connect summed.\hfill\break
HINT: If you have trouble see the proof of Theorem~\ref{thm:cstoob}.
\ehw

\begin{cor}[Giroux 2000, \cite{Giroux??}]
Let $a$ be any arc in $\Sigma,$ then 
\[
M_{S_{(\pm,a)}(\Sigma,\phi)} \text{ is diffeomorphic to } M_{(\Sigma,\phi)}
\]
and
\[
\xi_{S_{(+,a)}(\Sigma,\phi)} \text{ is isotopic to } \xi_{(\Sigma,\phi)}
\]
(where the corresponding manifolds are identified using the first diffeomorphism).
\end{cor}
\bbr
The contact structure $\xi_{S_{(-,a)}(\Sigma,\phi)}$ is not isotopic to $\xi_{(\Sigma,\phi)}$! One
can show that these contact structures are not even homotopic as plane fields.
\eer
\begin{proof}
The first statement was proven in the previous section. For the second statement recall if $(H^+,
\pi_+)$ is the open book for the positive Hopf link then $\xi_{(H^+,\pi_+)}$ is the standard contact
structure on $S^3.$ Thus
\[
\xi_{S_{(+,a)}(\Sigma,\phi)}=\xi_{(\Sigma,\phi)*(H^+, \pi_+)}=\xi_{(\Sigma,\phi)}\#\xi_{(H^+,\pi_+)}= \xi_{(\Sigma,\phi)}
\]
where all the equal signs mean isotopic. The last equality follows from the following exercise. 
\bhw
Show the contact manifold $(S^3,\xi_{std})$ is the union of two standard (Darboux) balls. 
\ehw
\bhw\label{ex:compBtight}
Show $\xi|_{M\setminus B}$ is tight. (If you do not know the definition of {\em tight} then see the beginning
of the next section.)\hfill\break
HINT: Try to show the contact structure pulled back to the cover $\Sigma\times\R$ of $M\setminus B$ is tight. This will be 
much easier after you know something about convex surfaces and, in particular, Giroux's tightness criterion which is 
discussed at the beginning of the next section.  
\ehw
\end{proof}

Note that we now have a well defined map 
\begin{equation}\label{eqn:psi}
\Psi:\mathcal{O} \to \mathcal{C}
\end{equation}
where 
\[
\mathcal{C}=\{\text{oriented contact structures on } M \text{ up to isotopy}\}
\] 
and 
\[
\mathcal{O}=\{\text{open book decompositions of } M \text{ up to positive stabilization}\}.
\]
In the next section we will show that $\Psi$ is onto and one-to-one.

\section{From contact structures to open books}\label{sec:proof}
We begin this section by recalling a few more basic facts and definitions from contact geometry.
Again for more details see \cite{Abetal, EtnyreCIntro, GeigesIntro}. This is not meant to be an 
introduction to contact geometry, but simply to remind the reader of some important facts or to
allow the reader with little background in contact geometry to follow some of the arguments below.
First some facts about Legendrian curves.
\begin{itemize}
\item A curve $\gamma$ in $(M,\xi)$ is Legendrian if it is always tangent to $\xi.$
\item Any curve may be $C^0$ approximated by a Legendrian curve.
\item If $\gamma\subset\Sigma$ is a simple closed Legendrian curve in $\Sigma$ then
$tw(\gamma,\Sigma)$ is the twisting of $\xi$ along $\gamma$ relative to $\Sigma.$ Said another way,
both $\xi $ and $\Sigma$ give $\gamma$ a framing (that is, a trivialization of its normal
bundle) by taking a vector field normal to $\gamma$ and tangent to $\xi$ or $\Sigma,$ respectively;
then $tw(\gamma,\Sigma)$ measures how many times the vector field corresponding to $\xi$ rotates as
$\gamma$ is traversed measured with respect to the vector field corresponding to $\Sigma.$
\end{itemize}

We now turn to surfaces in contact 3-manifolds and the fundamental dichotomy in 3-dimensional 
contact geometry: tight vs. overtwisted.
\begin{itemize}
\item If $\Sigma$ is a surface in $(M,\xi)$ then if at each point $x\in\Sigma$ we consider
$l_x=\xi_x\cap T_x\Sigma$ we get a singular line field on $\Sigma.$ This (actually any) line field
can be integrated to give a singular foliation on $\Sigma.$ This singular foliation is called the
\dfn{characteristic foliation} and is denoted $\Sigma_\xi.$
\item The contact structure $\xi$ is called \dfn{overtwisted} if there is an embedded disk $D$ such
that $D_\xi$ contains a closed leaf. Such a disk is called an overtwisted disk. If there are no
overtwisted disks in $\xi$ then the contact structure is called \dfn{tight}.
\item The standard contact structures on $S^3$ and $\R^3$ are tight.
\item If $\xi$ is a tight contact structure and $\Sigma$ is a surface with Legendrian boundary then
we have the weak-Bennequin inequality
\[
tw(\partial \Sigma,\Sigma)\leq \chi(\Sigma).
\]
\end{itemize}

We now begin a brief discussion of convex surfaces. These have proved to be an invaluable tool in
studying 3-dimensional contact manifolds.
\begin{itemize}
\item A surface $\Sigma$ in $(M,\xi)$ is called \dfn{convex} if there is a vector field $v$
transverse to $\Sigma$ whose flow preserves $\xi.$ A vector field whose flow preserves $\xi$ is
called a \dfn{contact vector field}.
\item Any closed surface is $C^\infty$ close to a convex surface. If $\Sigma$ has Legendrian
boundary such that $tw(\gamma,\Sigma)\leq 0$ for all components $\gamma$ of $\partial \Sigma$ then
after a $C^0$ perturbation of $\Sigma$ near the boundary (but fixing the boundary) $\Sigma$ will be
$C^\infty$ close to a convex surface.
\item Let $\Sigma$ be convex, with $v$ a transverse contact vector field. The set
\[
\Gamma_\Sigma=\{x\in \Sigma| v(x)\in \xi_x\}
\] 
is a multi-curve on $\Sigma$ and is called the
\dfn{dividing set}.
\item Let $\mathcal{F}$ be a singular foliation on $\Sigma$ and let $\Gamma$ be a multi-curve on
$\Sigma.$ The multi-curve $\Gamma$ is said to \dfn{divide} $\mathcal{F}$ if
\begin{enumerate}
\item $\Sigma\setminus\Gamma=\Sigma_+\coprod \Sigma_-$
\item $\Gamma$ is transverse to $\mathcal{F}$ and
\item there is a vector field $X$ and a volume form $\omega$ on $\Sigma$ so that
\begin{enumerate}
\item $X$ is tangent to $\mathcal{F}$ at non-singular points and $X=0$ at the singular points of
$\mathcal{F}$ (we summarize this by saying $X$ \dfn{directs} $\mathcal{F}$)
\item the flow of $X$ expands (contracts) $\omega$ on $\Sigma_+$ ($\Sigma_-$) and
\item $X$ points out of $\Sigma_+.$
\end{enumerate}
\end{enumerate}
\item If $\Sigma$ is convex then $\Gamma_\Sigma$ divides $\Sigma_\xi.$
\item On any compact subset of $\Sigma_+$ we can isotop $\xi$ to be arbitrarily close to $T\Sigma_+$
while keeping it transverse to $\Gamma_\Sigma.$
\item If $\Sigma$ is convex and $\mathcal{F}$ is any other foliation divided by $\Gamma_\Sigma$ then
there is a $C^0$ small isotopy, through convex surfaces, of $\Sigma$ to $\Sigma'$ so that
$\Sigma'_\xi=\mathcal{F}.$
\item If $\gamma$ is a properly embedded arc or a closed curve on $\Sigma,$ a convex surface, and
all components of $\Sigma\setminus \gamma$ contain some component of $\Gamma_\Sigma\setminus \gamma$
then $\Sigma$ may be isotoped through convex surfaces so that $\gamma$ is Legendrian. This is called
{\em Legendrian realization}.
\item If $\Sigma_1$ and $\Sigma_2$ are convex, $\partial \Sigma_1=\partial \Sigma_2$ is Legendrian
and the surfaces meet transversely, then the dividing curves interlace as shown in
Figure~\ref{fig:round} and we can round the corner to get a single smooth convex surface with
dividing curves shown in Figure~\ref{fig:round}.
\begin{figure}[ht]
  \relabelbox \small {\epsfysize=1.5in\centerline{\epsfbox{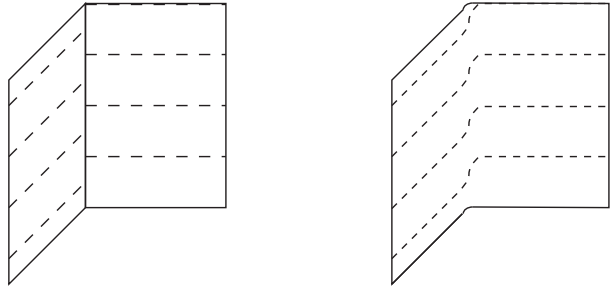}}} \endrelabelbox
        \caption{Two convex surfaces intersecting along their boundary, left, and the result of 
          rounding their corners, right.}
        \label{fig:round}
\end{figure}
\item If $\gamma$ is a Legendrian simple closed curve on a convex surface $\Sigma$ then
$tw(\gamma,\Sigma)=-\frac12 (\gamma\cap \Gamma_\Sigma).$
\item If $\Sigma$ is a convex surface then a small neighborhood of $\Sigma$ is tight if and only if
$\Sigma\not=S^2$ and no component of $\Gamma_\Sigma$ is a contractible circle or $\Sigma=S^2$ and
$\Gamma_\Sigma$ is connected. This is called Giroux's tightness criterion.
\end{itemize}
Now that we know about convex surfaces we can discuss a fourth way to say that a contact structure
is supported by an open book decomposition.
\begin{lem}
The contact structure $\xi$ on $M$ is supported by the open book decomposition $(B,\pi)$ if and only
if for every two pages of the open book that form a smooth surface $\Sigma'$ the contact structure
can be isotoped so that $\Sigma'$ is convex with dividing set $B\subset \Sigma'$ and $\xi$ is tight
when restricted to the components of $M\setminus \Sigma'.$
\end{lem}
This criterion for an open book to be supported by a contact structure is one of the easiest to
check and is quite useful in practice. This lemma is essentially due to Torisu \cite{Torisu00}.
\begin{proof}
Assume $\xi$ is supported by $(B,\pi).$ Let $V_0$ and $V_1$ be the closures of the complement of
$\Sigma'$ in $M.$
\bhw
Show $V_i$ is a handlebody.\hfill\break
HINT: Each $V_i$ is diffeomorphic to $\Sigma\times[0,1]$ where $\Sigma$ is a surface with boundary 
({\em i.e.} the page).
\ehw
It is not hard to show that $\Sigma'$ is convex. Indeed $\Sigma'$ is the union of two pages
$\Sigma_1$ and $\Sigma_2.$ Each $\Sigma_i$ has a transverse contact vector field $v_i.$ Along
$\partial \Sigma_i$ the $v_i$ point in opposite directions.
\bhw
Show $v_1$ and $-v_2$ can be altered in a neighborhood of $\partial \Sigma_1=\partial \Sigma_2=B$ so
that they give a contact vector field $v$ on $\Sigma'$ so that $B$ is the dividing set.\hfill\break
HINT: If you have trouble see \cite{Giroux91}.
\ehw 
In Exercise~\ref{ex:compBtight} you checked that $\xi$ restricted to $M\setminus B$ is tight. It is
easy to contact isotop $B$ to be disjoint from $V_i$ so $\xi$ restricted to $V_i$ is tight.

The other implication immediately follows from the next lemma.
\end{proof}
\begin{lem}[Torisu 2000, \cite{Torisu00}]
Given an open book decomposition $(B,\pi)$ of $M$ there is a unique contact structure $\xi$ that
makes $\Sigma'$ (the smooth union of two pages) convex with dividing set $B$ and that is tight when
restricted to each component of $M\setminus \Sigma'.$
\end{lem}
\begin{proof}[Sketch of Proof]
Let $\Sigma\subset \Sigma'$ be a page of the open book. Let $\alpha_1,\ldots, \alpha_n$ be a
collection of disjoint properly embedded arcs in $\Sigma$ that cut $\Sigma$ into a 2-disk. Since
each component of $M\setminus \Sigma'$ is a handlebody $V_i=\Sigma\times[0,1]$ we can consider the
disks $D_j=\alpha_j\times[0,1].$ These disks cut $V_i$ into a 3-ball. We can Legendrian realize
$\partial D_j$ on $\Sigma'$ and make all the disks $D_j$ convex. Now cutting $V_i$ along these disks
and rounding corners we get a tight contact structure on the 3-ball. Eliashberg \cite{Eliashberg92}
has shown that there is a unique tight contact structure on the 3-ball with fixed characteristic
foliation. From this it follows that there is a unique tight contact structure on $V_i$ with any fixed
characteristic foliation divided by $B.$ Finally this implies there is at most one contact structure
on $M$ satisfying the conditions in the lemma. The existence of one contact structure satisfying
these conditions is given by Theorem~\ref{twthm}.
\bhw
Fill in the details of this proof.\hfill\break
HINT: If you have trouble read the section on convex surfaces in \cite{EtnyreCIntro} or see 
\cite{Torisu00}.
\ehw
\end{proof}

We are now ready to show that the map $\Psi$ from open books to contact structures is onto.
\begin{thm}[Giroux 2000, \cite{Giroux??}]\label{thm:cstoob}
Every oriented contact structure on a closed oriented 3--manifold is supported by an open book 
decomposition.
\end{thm}
\begin{proof}
We begin the proof with a definition.
\begin{defn}
A \dfn{contact cell decomposition} of a contact 3--manifold $(M,\xi)$ is a finite CW-decomposition 
of $M$ such that
\begin{enumerate}
\item[(1)] the 1-skeleton is a Legendrian graph,
\item[(2)] each 2-cell $D$ satisfies $tw(\partial D, D)=-1,$ and
\item[(3)] $\xi$ is tight when restricted to each 3-cell.
\end{enumerate}
\end{defn}
\begin{lem}
Every closed contact 3-manifold $(M,\xi)$ has a contact cell decomposition.
\end{lem}
\begin{proof}
Cover $M$ by a finite number of Darboux balls (this is clearly possible since $M$ is compact). Note that
since Darboux balls are by definition contactomorphic to a ball in the standard contact structure on
$\R^3$ we know $\xi$ restricted to the Darboux balls is tight. Now take any finite CW-decomposition
of $M$ such that each 3-cell sits in some Darboux ball. Isotop the 1-skeleton to be Legendrian (this
can be done preserving the fact that 3-cells sit in Darboux balls). Note that we have a CW-decomposition
satisfying all but condition (2) of contact cell decomposition. To achieve this condition consider a
2-cell $D.$ By the weak-Bennequin inequality we have $tw(\partial D, D)\leq -1.$ Thus we can perturb
each 2-cell to be convex (care must be taken at the boundary of the 2-cells). Since $\Gamma_D$
contains no simple closed curves and $tw(\partial D, D)=-\frac12(\Gamma_D\cap \partial D)$ we know
that there are $\frac12 (\Gamma_D \cap \partial D)$ components to $\Gamma_D.$ If $tw(\partial D,
D)\not=-1$ there is more than one component to $\Gamma_D$ and we can thus use Legendrian realization
to realize arcs separating all the components of $\Gamma_D$ by Legendrian arcs. If we add these arcs
to the 2-skeleton then condition (2) of the definition is also satisfied.
\end{proof}

Suppose we have a contact cell decomposition of $(M,\xi).$ Denote its 1-skeleton by $G.$ Given the 
(or any) Legendrian graph $G$ the \dfn{ribbon of $G$} is a compact surface $R=R_G$ satisfying
\begin{enumerate}
\item $R$ retracts onto $G,$
\item $T_p R=\xi_p$ for all $p\in G,$
\item $T_p R\not = \xi_p$ for all $p\in R\setminus G.$
\end{enumerate}
Clearly any Legendrian graph has a ribbon. Let $B=\partial R$ and note that $B$ is a transverse link.
\begin{claim}
$B$ is the binding of an open book decomposition of $M$ that supports $\xi.$
\end{claim}
Clearly the proof follows form this claim. \end{proof}
\begin{proof}[Proof of Claim]
Since $B$ is a transverse link there is a contactomorphism from each component of a neighborhood
$N(B)$ of $B$ to an $\epsilon$-neighborhood of the $z$-axis in $(\R^3,\ker (dz+r^2\, d\theta))/\sim$
where $(r,\theta,z)\sim(r,\theta,z+1).$ Let $X(B)=\overline{M\setminus N(B)}$ be the complement of
$N(B)$ and $R_X=R\cap X(B).$ We can choose a neighborhood $N(R)=R_X\times[-\delta,\delta]$ of $R_X$
in $X(B)$ such that $\partial R_X\times\{pt\}$ thought of as sitting in $N(B)$ is a line with
constant $\theta$ value. Clearly $N(R)$ is an $R_X$ bundle over $[-\delta, \delta].$ Set
$X(R)=X(B)\setminus N(R).$ See Figure~\ref{fig:nbhd1}.
\begin{figure}[ht]
  \relabelbox \small {\epsfysize=1.7in\centerline{\epsfbox{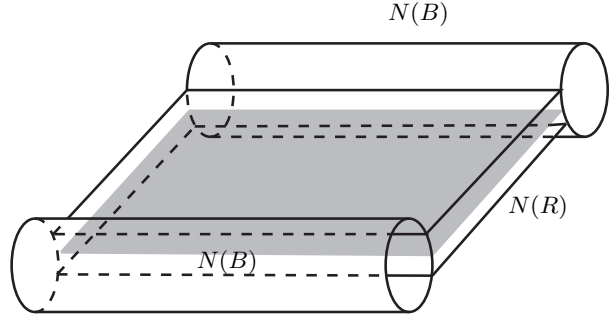}}} 
  \relabel{N1}{$N(B)$}
  \relabel{N2}{$N(B)$}
  \relabel{X}{$N(R)$}
  \endrelabelbox
        \caption{The neighborhoods $N(B)$ and $N(R)$. The grey part is $R_X.$}
        \label{fig:nbhd1}
\end{figure}

We first show that $X(R)$ is diffeomorphic to $R_X\times[0,1]$ and that $X(B)$ is formed by
identifying $R_X\times\{0\}$ ($R_X\times\{1\}$) in $X(R)$ and $R_X\times\{\delta\}$
($R_X\times\{-\delta\}$) in $N(R).$ Clearly this implies that $X(B),$ the complement of a
neighborhood of $B,$ is fibered and the fibration can be extended over $N(B)\setminus B$ so that the
boundary of the fibers is $B.$ Note that $\partial X(R)= A\cup F$ where $A=(\partial X(R))\cap N(B)$
is a disjoint union of annuli (one for each component of $N(B)$) that are naturally fibered by
circles of constant $\theta$ value in $N(B).$ The subsurface $F$ is defined to be the closure of the
complement of $A$ in $\partial X(R).$ Note that we can write $F=F^-\cup F^+,$ where $F^\pm$ is identified
with $R_X\times\{\pm\delta\}$ in $N(R).$
\bhw
Show that $\partial X(R)$ is a convex surface with dividing set $\Gamma_{\partial X(R)}$ equal to the
union of the cores of $A$ and such that $F^\pm\subset (\partial X(R))_\pm.$ (Note, $\partial X(R)$
is only piecewise smooth, but if we rounded the edges it would be convex.)
\ehw
\bbr
Throughout this part of the proof we will be discussing manifolds whose boundaries have corners. We
do not want to smooth the corners. However, sometimes to understand the annuli $A$ better we will
think about rounding the corners, but once we have understood $A$ sufficiently we actually will not
round the corners.
\eer
Let $D_1,\ldots, D_k$ be the two cells in the contact cell decomposition of $(M,\xi).$ Recall that
$\partial D_i$ is Legendrian and has twisting number $-1.$ Thus since $R$ twists with the contact
structure along the 1-skeleton $G$ we can assume that $B$ intersects $D_i$ exactly twice for all
$i.$ Let $D_i'=D_i\cap X(R).$
\bhw
Show that it can be arranged that the $D_i'$'s intersect the region $A\subset \partial X(R)$ in exactly
two properly embedded arcs, and each arc runs from one boundary component of $A$ to another.
\ehw
\bhw
Show that the interior of $X(R)$ cut along all the $D_i'$'s is homeomorphic to $M$ minus its 2-skeleton.
That is, $X(R)$ cut along the $D_i'$'s is a union of balls.  
\ehw
Using the Legendrian realization principle we can assume $\partial D_i'$ is Legendrian. (Again as in
the remark above we only Legendrian realize $\partial D_i'$ to see what happens to the dividing
curves when we cut $X(R)$ along these disks. Once we have seen this we don't actually do the
Legendrian realization.) Let's consider what happens to $X(R)$ when we cut along $D_1';$ denote the
resulting manifold by $X_1.$ Note: in $\partial X_1$ there are two copies of $D_1'.$ Let $A_1$ be $A$
sitting in $\partial X_1$ union the two copies of $D_1.$ See Figure~\ref{fig:cutx}.
\begin{figure}[ht]
  \relabelbox \small {\epsfysize=5in\centerline{\epsfbox{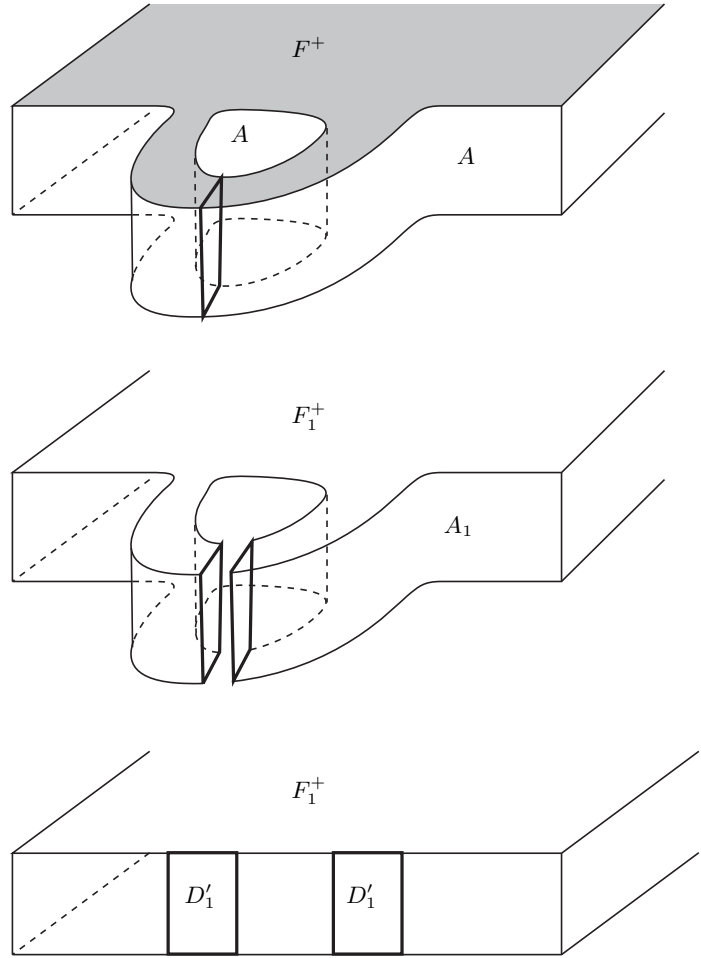}}} 
  \relabel{A1}{$A$}
  \relabel{A2}{$A$}
  \relabel{A}{$A_1$}
  \relabel{Fp}{$F^+$}
  \relabel{Fpp}{$F^+_1$}
  \relabel{Fppp}{$F^+_1$}
  \relabel{D1}{$D'_1$}
  \relabel{D1p}{$D'_1$}
  \endrelabelbox
        \caption{Top picture is $X(R)$ with the boundary of $D_1'$ drawn darkly. The middle picture 
          is $X_1$ right after the cut and the bottom picture is $X_1$ after isotoping a little.}
        \label{fig:cutx}
\end{figure}
\bhw
Show that the components of $A_1$ are annuli and that they have a natural fibration by $S^1$ that is naturally 
related to the fibration on $A.$
\ehw 
Note that $\partial X_1\setminus A_1$ naturally breaks into two surfaces $F_1^+$ and $F_1^-,$ where
$F_1^\pm$ is obtained from $F^\pm$ by cutting along a properly embedded arc.
\bhw
Show $\partial X_1$ is convex (once the corners are rounded) and its dividing set is the union of
the cores of $A_1$ and $F_1^\pm\subset (\partial X_1)_\pm.$
\ehw
If we continue to cut along the $D_i'$'s we eventually get to $X_k$ once we have cut along all the
disks. From above we know $X_k$ is a disjoint union of balls (all contained in 3-cells of our
contact cell decomposition). Moreover, on $\partial X_k$ we have that $A_k$ a union of annuli whose cores
give the dividing curves for $\partial X_k.$ By the definition of contact cell decomposition we know that
the contact structure when restricted to each component of $X_k$ is tight. Thus we know $A_k$ has
exactly one component on the boundary of each component of $X_k.$ Thus each component of $X_k$ is a
ball $B^3$ with an annulus $S$ that has a natural fibration by circles. Clearly $B^3$ has a natural
fibration by $D^2$'s that extends the fibration of $S$ by circles. That is, $B^3=D^2\times[0,1]$ with
$(\partial D^2)\times[0,1]=S.$
\bhw
Show that as we glue $X_k$ together along the two components of $D_k'$ in $\partial X_k$ to get $X_{k-1}$
we can glue the fibration of $X_k$ by $D^2$'s together to get a fibration of $X_{k-1}$ by surfaces
that extend the fibration of $A_{k-1}$ by circles.
\ehw
Thus, continuing in this fashion we get back to $X(R)$ and see that it is fibered by surfaces that
extend the fibration of $A$ by circles. This clearly implies that $X(R)=F^-\times[0,1]=R_X\times
[0,1]$ and the surfaces $R_X\times\{0\}$ and $R_X\times\{1\}$ are glued to the boundary of $N(R)$ as
required above. Hence we have shown that $X(B)$ is fibered over the circle by surfaces diffeomorphic
to $R$ and that the fibers all have boundary $B.$ That is, we have demonstrated that $B$ is the binding of
an open book.

We now must show that this open book supports the contact structure $\xi.$ Looking back through the
proof it is not hard to believe that one may isotop the contact planes to be arbitrarily close to
the pages of the open book, but it seems a little difficult to prove this directly. We will show the
open book is compatible with the contact structure by showing that there is a Reeb vector field that
is tangent to the binding and transverse to the pages. Recall that the neighborhood $N(B)$ of the
binding is contactomorphic to an $\epsilon$-neighborhood of the $z$-axis in $(\R^3,\ker (dz+r^2\,
d\theta))/\sim$ where $(r,\theta,z)\sim(r,\theta,z+1).$ Moreover, we can assume the pages intersect 
this neighborhood as the constant $\theta$ annuli.
\bhw
In the explicit model for $N(B)$ find a Reeb vector field that is tangent to the binding and
positively transverse to the pages of the open book in the neighborhood. Also make sure the boundary
of $N(B)$ is preserved by the flow of the Reeb field.
\ehw
We can think of the Reeb fields just constructed as giving a contact field in the neighborhood of
the boundary of $R_X$ (recall this is the ribbon of the Legendrian 1-skeleton $G$ intersected with
the complement of $N(B)$).
\bhw
Show that this contact vector field defined in a neighborhood of $\partial R_X$ can be extended to a
contact vector field $v$ over the rest of $R_X$ so that it is transverse to $R_X$ and there are no
dividing curves. (This is OK since $R_X$ is not a closed surface.) Note that since there are no
dividing curves $v$ is also transverse to $\xi.$
\ehw
Use $v$ to create the neighborhood $N(R)$ of $R_X.$ Since $v$ is never tangent to the contact planes 
along $R_X$ we can assume that this is the case in all of $N(R).$ 
\bhw
Show that a contact vector field which is never tangent to the contact planes is a Reeb vector field. 
\ehw
Thus we have a Reeb vector field defined on $N(B)\cup N(R)$ that has the desired properties.

We now need to extend the Reeb vector field $v$ over $X(R).$ From the construction of $v$ we can
assume we have $v$ defined near the boundary of $X(R)$ and as a vector field defined there it
satisfies the following:
\begin{enumerate}
\item $v$ is tangent to $A\subset \partial X(R).$
\item There is a neighborhood $N(A)$ of $A$ in $X(R)$ such that $\partial N(A)= A\cup A^+\cup
A^-\cup A'$ where $A^\pm=N(A)\cap F^\pm$ and $A'$ is a parallel copy of $A$ on the interior of
$X(R).$ $v$ is defined in $N(A),$ tangent to $A\cup A',$ $\pm v$ points transversely out of $N(A)$
at $A^\pm$ and $v$ is transverse to the pages of the open book intersected with $N(A).$ Moreover the
flow of $v$ will take $A^-$ to $A^+.$
\item $\pm v$ points transversely out of $X(R)$ along $F^\pm.$
\end{enumerate}

We now want to construct a model situation into which we can embed $X(R).$ To this end let
$\Sigma=R_X\cup A' \cup -R_X,$ where $A'=(\partial R_X)\times [0,1]$ and the pieces are glued
together so that $\Sigma$ is diffeomorphic to the double of $R_X.$ On $\Sigma$ let $\mathcal{F}$ be
the singular foliation $(R_X)_\xi$ on each of $R_X$ and $-R_X$ and extend this foliation across $A'$
so that it is non singular there and the leaves of the foliation run from one boundary component to
another. Let $\Gamma$ be the union of the cores of the annuli that make up $A'.$ It is easy to see
that $\mathcal{F}$ is divided by $\Gamma.$ Given this one can create a vertically invariant contact
structure $\xi'$ on $\Sigma\times\R$ such that $(\Sigma\times\{t\})_{\xi'}=\mathcal{F}$ and the
dividing set on $\Sigma\times\{t\}$ is $\Gamma,$ for all $t\in \R.$ (See \cite{Giroux91}.) Note that
$\frac{\partial}{\partial t}$ restricted to $R_X\times\R$ is a Reeb vector field since it is a
contact vector field and positively transverse to $\xi'$ in this region. Pick a diffeomorphism
$f:F^-\to (R_X\times\{0\}\subset \Sigma\times\{0\})$ that sends $(R_X)_\xi$ to $\mathcal{F}.$
(Recall $F^+\cup F^-=(\partial X(R))\setminus A.$)
\bhw
Show that the flow of $v$ on $R_X$ and $\frac{\partial}{\partial t}$ on $\Sigma\times\R$ allow you
to extend $f$ to a contact embedding of $N(A)$ into $\Sigma\times\R.$
\ehw
Thus we can use the flow of $v$ and $\frac{\partial}{\partial t}$ to extend $f$ to a contact
embedding of a neighborhood of $\partial X(R)$ in $X(R)$ into $\Sigma\times\R.$
\bhw
Make sure you understand how to get the embedding near $F^+.$\hfill\break
HINT: From the previous exercise we have a neighborhood of the boundary of $F^+$ embedded into
$R_X\times\{t_0\},$ for some $t_0.$ Show that there is an obvious way to extend this to an embedding of
all of $F^+$ to $R_X\times\{t_0\}.$
\ehw
Of course this extension of $f,$ which we also call $f,$ takes the Reeb field $v$ to the Reeb field
$\frac{\partial}{\partial t}.$ We can clearly extend $f$ to an embedding, but not necessarily a
contact embedding, of all of $X(R)$ into $\Sigma\times\R.$ The following exercises allow us to
isotop $f,$ relative to a neighborhood of the boundary, to a contact embedding and thus we may
extend $v$ to all of $X(R)$ by $\frac{\partial}{\partial t}.$ This gives us a Reeb vector field on
$M$ which demonstrates that the open book supports $\xi.$
\bhw
Let $H$ be a handlebody and $D_1,\ldots, D_g$ be properly embedded disks that cut $H$ into a 3-ball.
Given any singular foliation $\mathcal{F}$ on the boundary of $H$ that is divided by $\Gamma$ for
which $\partial D_i\cap \Gamma=2,$ for all $i,$ then there is at most one tight contact structure on
$H,$ up to isotopy, that induces $\mathcal{F}$ on $\partial H.$\hfill\break
HINT: This is a simple exercise in convex surface theory. See \cite{EtnyreCIntro}.
\ehw
\bhw
Show the contact structure $\xi'$ on $\Sigma\times \R$ is tight. (This is easy using Giroux's
tightness criterion.) Also show the contact structure $\xi$ restricted to $X(R)$ is tight.
\hfill\break 
HINT: The second part is not so easy. The idea is that if you can cut up a handlebody by disks, as in
the previous exercise, and the 3-ball you end up with has a tight contact structure on it, then the
original contact structure on the handlebody is tight. See \cite{Honda02}.
\ehw
\end{proof}

\begin{figure}[ht]
  \relabelbox \small {\epsfysize=2.2in\centerline{\epsfbox{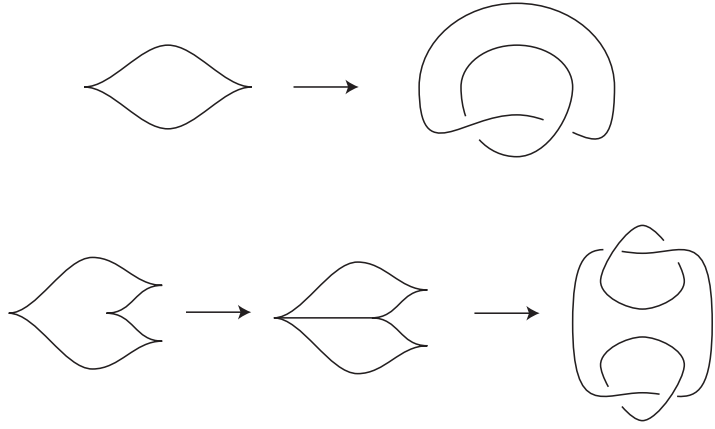}}} 
  \endrelabelbox
        \caption{On the top left we have a Legendrian unknot that is the 1-skeleton of a contact 
          cell decomposition of $S^3.$ The resulting open book is shown on the upper right. On the bottom 
          left we start with a Legendrian unknot, moving to the right we add a Legendrian arc to 
          get the 1-skeleton of a contact cell decomposition. The bottom right shows the resulting open 
          book.}
        \label{fig:rx}
\end{figure}

We have the following immediate useful corollaries.
\begin{cor}\label{Lonpage}
If $L$ is a Legendrian link in $(M,\xi)$ then there is an open book decomposition supporting $\xi$
such that $L$ sits on a page of the open book and the framing given by the page and by $\xi$ agree.
\end{cor}
\begin{proof}
Simply include the Legendrian link $L$ in the 1-skeleton of the contact cell decomposition.
\end{proof}
\bex
Figure~\ref{fig:rx} illustrates Corollary~\ref{Lonpage} for two knots in $S^3$ with its standard
contact structure.
\eex

Using positive stabilizations we can see the following.
\begin{cor}
Any contact manifold is supported by an open book with connected binding.
\end{cor}
\begin{cor}[Contact Bing]
A contact manifold $(M,\xi)$ is the standard tight contact structure on $S^3$ if and only if every
simple closed curve is contained in a Darboux ball.
\end{cor}
\bhw
Prove this last corollary.\hfill\break
HINT: There is a unique tight contact structure on $B^3$ inducing a fixed characteristic foliation 
on the boundary \cite{Eliashberg92}.
\ehw

\begin{thm}[Giroux 2002, \cite{Giroux??}]\label{thm:stab}
Two open books supporting the same contact manifold $(M,\xi)$ are related by positive stabilizations.
\end{thm}
To prove this theorem we need the following lemma.
\begin{lem}\label{makecell}
Any open book supporting $(M,\xi),$ after possibly positively stabilizing, comes from a contact cell
decomposition.
\end{lem}
\begin{proof}
Let $\Sigma$ be a page of the open book and $G$ be the core of $\Sigma.$ That is, $G$ is a graph
embedded in $\Sigma$ onto which $\Sigma$ retracts. We can Legendrian realize $G.$ 
\bbr\label{extralerp}
The Legendrian realization principle is for curves, or graphs, on a closed convex surface or a
convex surface with Legendrian boundary. The pages of an open book are convex but their boundary is
transverse to the contact structure so we cannot apply the Legendrian realization principle as it is
usually stated. Nonetheless since we can keep the characteristic foliation near the boundary fixed
while trying to realize a simple closed curve or graph, we can still realize it. But recall the curve
or graph must be non-isolating. In this context this means that all components of the complement of
the curve in the surface should contain a boundary component. To see this review Giroux's proof of
realization.
\eer
Note that $\Sigma$ is the ribbon of $G.$ Let $N$ be a neighborhood of $\Sigma$ such that $\partial
\Sigma\subset \partial N.$ Let $\alpha_i$ be a collection of properly embedded arcs on $\Sigma$ that
cut $\Sigma$ into a disk. Let $\widetilde{A}_i$ be $\alpha_i\times[0,1]$ in $M\setminus
\Sigma=\Sigma \times[0,1]$ and $A'_i=\widetilde{A}_i\cap \overline{(M\setminus N)}.$ Note that $A'_i$
intersects $\partial \Sigma $ on $\partial N$ exactly twice. Thus if we extend the $A_i$'s into $N$
so their boundaries lie on $G$ then the twisting of $\Sigma,$ and hence $\xi,$ along $\partial A_i$
with respect to $A_i$ is $-1, 0$ or 1. If all the twisting is $-1$ then we have a contact cell
decomposition (recall that the contact structure restricted to the complement of $\Sigma$ is tight). Thus
we just need to see how to reduce the twisting of $\xi$ along $\partial A_i.$

Suppose $\partial A_i$ has twisting 0. Positively stabilize $\Sigma$ as shown in
Figure~\ref{fig:cfram}.
\begin{figure}[ht]
  \relabelbox \small {\epsfysize=1.4in\centerline{\epsfbox{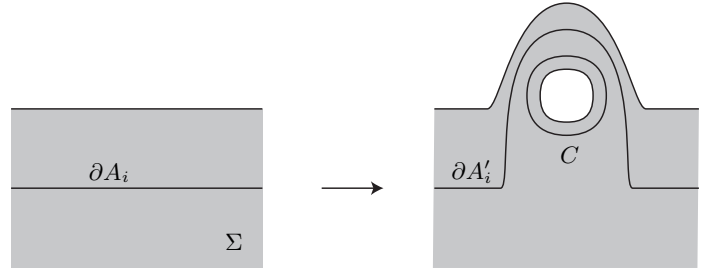}}} 
  \relabel{A}{$\partial A_i$}
  \relabel{N}{$\partial A_i'$}
  \relabel{S}{$\Sigma$}
  \relabel{C}{$C$}
  \endrelabelbox
        \caption{On the left is part of $\Sigma$ and $\partial A_i$ near $\partial \Sigma.$ On the 
          right is the stabilized $\Sigma, C$ and $\partial A_i'.$ (This picture is abstract. If we 
          were drawing it ambiently the added 1-handle would have a full left-handed twist in it.)}
        \label{fig:cfram}
\end{figure}
Note the curve $C$ shown in the picture can be assumed to be Legendrian and bounds a disk $D$ in
$M.$ Now isotop $G$ across $D$ to get a new Legendrian graph with all the $A_j$'s unchanged except
that $A_i$ is replace with the disk $A_i'$ obtained from $A_i$ by isotoping across $D.$ We also
add $C$ to $G$ and add $D$ to the 2-skeleton.
\bhw
Show that the twisting of $\xi$ along $\partial A_i'$ is one less than the twisting along $\partial A_i.$
\ehw
Clearly the twisting of $\xi$ along $D$ is $-1.$ Thus we can reduce the twisting of $\xi$ along
$\partial A_i$ as needed and after sufficiently many positive stabilizations we have an open book
that comes from a contact cell decomposition.
\end{proof}

\begin{proof}[Proof of Theorem~\ref{thm:stab}]
Given two open books $(B,\pi)$ and $(B',\pi')$ supporting $(M,\xi)$ we can assume they both come
from contact cell decompositions by using Lemma~\ref{makecell}. Now, given two contact cell
decompositions one can show that they are related by a sequence of the following:
\begin{enumerate}
\item A subdivision of a 2-cell by a Legendrian arc intersecting the dividing set one time.
\item Add a 1-cell $c'$ and a 2-cell $D$ so that $\partial D=c'\cup c$ where $c$ is part of the
  original 1-skeleton and $tw(\partial D, D)=-1.$
\item Add a 2-cell $D$ whose boundary is already in the 1-skeleton and $tw(\partial D, D)=-1.$
\end{enumerate}

Thus the theorem follows from the following exercises.
\bhw
Show that (3) does not change the open book associated to the cell decomposition.
\ehw
\bhw
Show that (1) and (2) positively stabilize the open book associated to the cell decomposition.
\hfill\break
HINT: Show that an arc $a$ is added to the 1-skeleton and a disk $D$ to the 2-skeleton so that $\partial
D$ is part of the old 1-skeleton union $a,$ $\Gamma_D$ is a single arc and $a\cap \Gamma_D$ is one
point. Show that adding such an arc to the one skeleton is equivalent to positively stabilizing the open
book.
\ehw
\end{proof}

\section{Symplectic cobordisms and caps.}\label{sec:apps}
A contact manifold $(M,\xi)$ is called \dfn{(weakly) symplectically fillable} if there is a compact
symplectic 4--manifold $(X,\omega)$ such that $\partial X=M$ and $\omega|_\xi\not=0.$ Many
applications of contact geometry to topology (see \cite{KronheimerMrowka, OzsvathSzabo} and the
discussion in the introduction) rely on the following theorem.
\begin{thm}[Eliashberg 2004 \cite{Eliashberg04}; Etnyre 2004 \cite{Etnyre04}]\label{thm:caps}
If $(X,\omega)$ is a symplectic filling of $(M,\xi)$ then there is a closed symplectic manifold
$(W,\omega')$ and a symplectic embedding $(X,\omega)\to (W,\omega').$
\end{thm}
Partial results aimed towards this theorem were obtained by many people. Specifically, Lisca and
Mati\'c established this result for Stein fillable manifolds in \cite{LiscaMatic97} and later work of
Akbulut and Ozbagci \cite{AO} coupled with work of Plamenevskaya \cite{Plam} provided an alternate
proof in this case (for unfamiliar terminology see the next paragraph). For strongly fillable
manifolds this was proven by Gay in \cite{Gay02} and follows trivially from Theorem 1.3 in
\cite{EtnyreHonda02a}. The full version of this theorem also follows fairly easily from
\cite{Stipsicz03}.

In the process of proving Theorem~\ref{thm:caps} we will need to take a few detours. The first
concerns various types of symplectic fillings and the second concerns Legendrian/contact surgery.
These two detours occupy the next two subsections. We return to the proof of Theorem~\ref{thm:caps}
in Subsection~\ref{subsec:capspf}. In the final section we discuss the relation between open book
decompositions and overtwistedness.

\subsection{Symplectic fillings}\label{subsec:sfill}
A contact manifold $(M,\xi)$ is said to be \dfn{strongly symplectically filled} by the symplectic
manifold $(X,\omega)$ if $X$ is compact, $\partial X=M$ and there is a vector field $v$ transversely
pointing out of $X$ along $M$ such that the flow of $v$ dilates $\omega$ (that is to say the Lie
derivative of $\omega$ along $v$ is a positive multiple of $\omega$). The symplectic manifold
$(X,\omega)$ is said to have convex boundary if there is a contact structure $\xi$ on $\partial X$
that is strongly filled by $(X,\omega).$ We say that $(X,\omega)$ is a \dfn{strong concave filling}
if $(X,\omega)$ and $v$ are as above except that $v$ points into $X.$ Note that given a symplectic
manifold $(X,\omega)$ with a dilating vector field $v$ transverse to its boundary then
$\iota_v\omega$ is a contact form on $\partial X.$ If $v$ points out of $X$ then the contact form
gives an oriented contact structure on $\partial X$ and if $v$ points into $X$ then it gives an
oriented contact structure on $-\partial X.$ (Recall Remark~\ref{rem:poscont})

Given a contact manifold $(M,\xi)$ let $\alpha$ be a contact form for $\xi,$ then consider $W=M\times
\R$ and set $\omega_W=d(e^t\alpha),$ where $t$ is the coordinate on $\R.$ It is easy to see that
$\omega_W$ is a symplectic form on $W$ and the vector field $v=\frac{\partial}{\partial t}$ is a
dilating vector field for $\omega_W.$ The symplectic manifold $(W,\omega_W)$ is called the 
\dfn{symplectization} of $(M,\xi).$
\bhw
Given any other contact form $\alpha'$ for $\xi$ (note that this implies that $\alpha'=g\alpha$ for some
function $g:M\to (0,\infty)$) show there is some function $f$ such that $\alpha'=F^*(\iota_v
\omega_W)$ where $F:M\to W:x\mapsto (x,f(x)).$
\ehw

\begin{figure}[ht]
  \relabelbox \small {\epsfysize=2in\centerline{\epsfbox{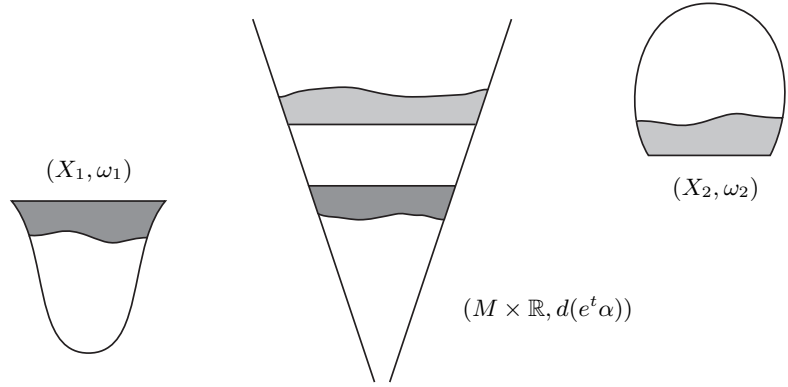}}} 
  \relabel{S}{$(M\times\R, d(e^t \alpha))$}
  \relabel{CC}{$(X_2,\omega_2)$}
  \relabel{CV}{$(X_1,\omega_1)$}
  \endrelabelbox
       \caption{The symplectization of $(M,\xi),$ middle, and a symplectic manifold with convex, 
         left, and concave, right, boundary.}
       \label{fig:symp}
\end{figure}

It can be shown that if $(X,\omega)$ is a strong symplectic filling (strong concave filling) of
$(M,\xi)$ then there is a neighborhood $N$ of $M$ in $X,$ a function $f,$ a one-sided neighborhood
$N_f$ of the graph of $f$ in $W$ with $N_W$ lying below (above) the graph and a symplectomorphism
$\psi:N_W\to N.$ See Figure~\ref{fig:symp}.
Thus we have a model for a neighborhood of a contact manifold in a strong symplectic filling.
\bhw
If $(X_1,\omega_1)$ is a strong symplectic filling of $(M,\xi)$ and $(X_2,\omega_2)$ is a strong
concave filling of $(M,\xi)$ then show $X=X_1\cup X_2$ has a symplectic structure $\omega$ such that
$\omega|_{X_1}=\omega_1$ and $\omega|_{X_2\setminus N}=c\omega_2$ where $N$ is a neighborhood of $\partial X_2$
in $X_2$ and $c>0$ is a
constant.\hfill\break
HINT: Look at Figure~\ref{fig:symp}.
\ehw
Thus we can use strong symplectic fillings to glue symplectic manifolds together. This is not, in 
general, possible with a weak symplectic filling.

Recall that a Stein manifold is a triple $(X,J,\psi)$ where $J$ is a complex structure on $X$ and
$\omega_\psi(v,w)=-d(d\psi \circ J)(v,w)$ is non-degenerate. A contact manifold $(M,\xi)$ is called
\dfn{Stein fillable} (or \dfn{holomorphically fillable}) if there is a Stein manifold $(X,J,\psi)$
such that $\psi$ is bounded from below, $M$ is a non-critical level of $\psi$ and $-(d\psi\circ J)$
is a contact form for $\xi.$ It is customary to think of $X$ as $\psi^{-1}((-\infty, c])$ where
$M=\psi^{-1}(c).$ Thus we can think of $X$ as a compact manifold (Stein manifolds themselves are
never compact).

In \cite{Eliashberg90, Weinstein91} it was shown how to attach a 1-handle to the boundary of a
symplectic manifold with convex boundary and extend the symplectic structure over the 1-handle so as
to get a new symplectic manifold with convex boundary. They also showed the same could be done when
a 2-handle is attached along a Legendrian knot with framing one less than the contact framing. In
fact we have the following characterization of Stein manifolds.
\begin{thm}[Eliashberg 1990, \cite{Eliashberg90}]\label{stein}
A 4--manifold $X$ is Stein if and only if $X$ has a handle decomposition with only 0-handles,
1-handles and 2-handles attached along Legendrian knots with framing one less than the contact
framing. \qed
\end{thm}

Summarizing the relations between various notions of filling and tightness we have
$$
\mbox{Tight}\quad \supset \quad \mbox{Weakly Fillable}\quad \supset \quad\mbox{Strongly Fillable} 
\quad \supset \quad\mbox{Stein Fillable}.
$$
The first two inclusions are strict, see \cite{EtnyreHonda02a} and \cite{Eliashberg96} respectively.
It is unknown whether or not the last inclusion is strict. We have the following useful fact.
\begin{thm}[Eliashberg 1991, \cite{Eliashberg91}; Ohta and Ono 1999, \cite{OO}]\label{filleqhom}
If $M$ is a rational homology sphere then any weak filling of $(M,\xi)$ can be deformed into a
strong filling.\qed
\end{thm}

\subsection{Contact surgery}
Let $L$ be a Legendrian knot in a contact 3--manifold $(M,\xi).$ It is well known (see
\cite{EtnyreCIntro, GeigesIntro}) that $L$ has a neighborhood $N_L$ that is contactomorphic to a
neighborhood of the $x$-axis in
\[
(\R^3,\ker(dz -y dx))/\sim,
\]
where $\sim$ identifies $(x,y,z)$ with $(x+1,y,z).$ With respect to these coordinates on $N_L$ we
can remove $N_L$ from $M$ and topologically glue it back with a $\pm 1$-twist (that is, we are doing
Dehn surgery along $L$ with framing the contact framing $\pm 1$). Call the resulting manifold
$M_{(L,\pm 1)}.$ There is a unique way, up to isotopy, to extend $\xi|_{M\setminus N_L}$ to a
contact structure $\xi_{(L,\pm 1)}$ over all of $M_{(L,\pm 1)}$ so that $\xi_{(L,\pm 1)}|_{N_L}$ is
tight (see \cite{Honda00}). The contact manifold $(M,\xi)_{(L,\pm 1)}=(M_{(L,\pm 1)},\xi_{(L,\pm
1)})$ is said to be obtained from $(M,\xi)$ by \dfn{$\pm 1$-contact surgery} along $L.$ It is
customary to refer to $-1$-contact surgery along $L$ as \dfn{Legendrian surgery} along $L.$
\begin{quest}
Is $\xi_{(L,-1)}$ tight if $\xi$ is tight?
\end{quest}
If the original contact manifold $(M,\xi)$ is not closed then it is known that the answer is
sometimes NO, see \cite{Honda02}. But there is no known such example on a closed manifold. It is
known, by a combination of Theorems~\ref{stein} and \ref{filleqhom}, that Legendrian surgery (but
not $+1$-contact surgery!) preserves any type of 
symplectic fillability. (Similarly, $+1$-contact surgery preserves non-fillability.) We have the following
result along those lines.
\begin{thm}[Eliashberg 1990, \cite{Eliashberg90}; Weinstein 1991, \cite{Weinstein91}]
\label{surgeryandcobord}
Given a contact 3--manifold $(M,\xi)$ let $(W=M\times [0,1], \omega=d(e^t \alpha))$ be a piece of
the symplectization of $(M,\xi)$ discussed in the last section. Let $L$ be a Legendrian knot sitting
in $(M,\xi)$ thought of as $M\times\{1\}.$ Let $W'$ be obtained from $W$ by attaching a 2-handle to
$W$ along $L\subset M\times\{1\}$ with framing one less than the contact framing. Then $\omega$ may
be extended over $W'$ so that the upper boundary is still convex and the induced contact manifold is
$(M_{(L,-1)},\xi_{(L,-1)}).$ Moreover, if the 2-handle was added to a Stein filling (respectively
weak filling, strong filling) of $(M,\xi)$ then the resulting manifold would be a Stein filling
(respectively weak filling, strong filling) of $(M_{(L,-1)},\xi_{(L,-1)}).$\qed
\end{thm}

We now want to see how contact surgery relates to open book decompositions. The main result along
these lines is the following.
\begin{thm}\label{surgeryonbook}
Let $(\Sigma,\phi)$ be an open book supporting the contact manifold $(M,\xi).$ If $L$ is a
Legendrian knot on the page of the open book then
\[
(M,\xi)_{(L,\pm 1)}= (M_{(\Sigma, \phi\circ D_L^{\mp})}, \xi_{(\Sigma, \phi\circ D_L^{\mp})}).
\]
\end{thm}
\begin{proof}
We begin by ignoring the contact structures and just concentrating on the manifold. We have a simple
closed curve $L$ on the page $\Sigma$ of the open book. Recall $M\setminus \text{nbhd }B$ is the
mapping cylinder $\Sigma_\phi.$ We will think of $L$ as sitting on $\Sigma\times\{\frac12\}$ in
$\Sigma\times[0,1],$ then by moding out by the identification $(\phi(x),0)\sim(x,1)$ we will have $L$
on a page in $\Sigma_\phi\subset M.$ 
\bhw
Show that cutting $\Sigma_\phi$ open along $\Sigma\times\{\frac12\}$ and regluing using $D_L^\pm$
will give you a manifold diffeomorphic to $\Sigma_{\phi\circ D^\pm_L}.$
\ehw

Let $N_\Sigma$ be a closed tubular neighborhood of $L$ in $\Sigma\times\{\frac12\}.$ Then a
neighborhood $N$ of $L$ in $M$ is given by $N_\Sigma\times [\frac12-\epsilon, \frac12+\epsilon].$ We
can assume the support of $D^\mp_L$ is in $N_\Sigma.$ Thus if $N'$ is neighborhood of $L$ in
$\Sigma_{\phi\circ D^\mp_L}$ corresponding to $N$ in $\Sigma_\phi$ then
\[
\Sigma_\phi\setminus N = \Sigma_{\phi\circ D^\mp_L}\setminus N'.
\]
So clearly $M_{(\Sigma,\phi)}\setminus N$ is diffeomorphic to $M_{(\Sigma,\phi\circ
D^\mp_L)}\setminus N'.$ Said another way $M_{(\Sigma,\phi\circ D^\mp_L)}$ is obtained from
$M_{(\Sigma,\phi)}$ by removing a solid torus and gluing it back in, {\em i.e.} by a Dehn surgery
along $L.$ We are left to see that the Dehn surgery is a $\pm 1$ Dehn surgery with respect to the
framing on $L$ coming from the page on which it sits. See Figure~\ref{fig:DequalD}
\begin{figure}[ht]
  \relabelbox \small {\epsfysize=2in\centerline{\epsfbox{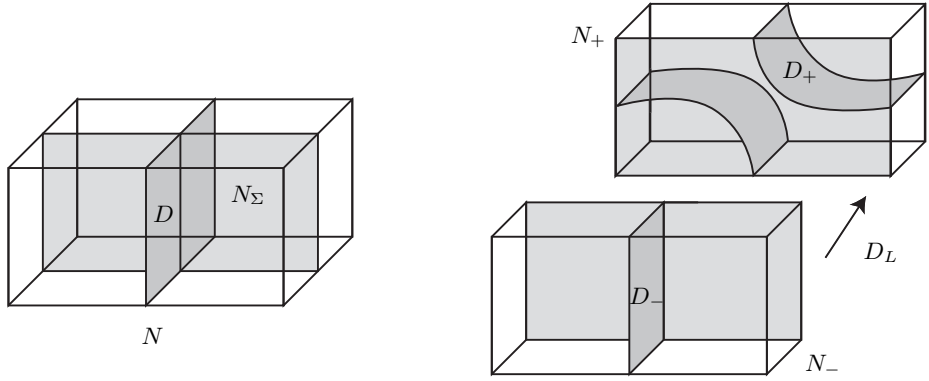}}} 
  \relabel{D}{$D$}
  \relabel{NS}{$N_\Sigma$}
  \relabel{N}{$N$}
  \relabel{L}{$D_L$}
  \relabel{DP}{$D_+$}
  \relabel{DM}{$D_-$}
  \relabel{NP}{$N_+$}
  \relabel{M}{$N_-$}
  \endrelabelbox
        \caption{On the right is the neighborhood $N$ with its meridional disk $D.$ On the right is 
        the neighborhood $N'.$ (The right- and left-hand sides of each cube are identified to get a
        solid torus.)}
        \label{fig:DequalD}
\end{figure}
while reading the rest of this paragraph. Note that we get $N'$ from $N$ by cutting $N$ along
$N_\Sigma$ and regluing using $D_L^\mp.$ We get a meridional disk $D$ for $N$ by taking an arc $a$
on $N_\Sigma\times\{\frac12-\epsilon\}$ running from one boundary component to the other and setting
$D=a\times[\frac12-\epsilon, \frac12+\epsilon].$ Let $N_-=N_\Sigma\times[\frac12-\epsilon, \frac12]$
and $N_+=N_\Sigma\times[\frac12,\frac12+\epsilon].$ Thus we get $N'$ by gluing
$\Sigma\times\{\frac12\}$ in $N_+$ to $\Sigma\times\{\frac12\}$ in $N_-.$ Set $D_-=D\cap N_-.$ Then
$a_-=D_-\cap \Sigma\times\{\frac12\}$ in $N_-$ is taken to $a_+\subset
\Sigma\times\{\frac12\}\subset N_+.$ So $D',$ the meridional disk in $N',$ is $D_-\cup
a_+\times[\frac12-\epsilon,\frac12].$ Thus on $\partial M_{(\Sigma,\phi)}\setminus N=\partial
M_{(\Sigma,\phi\circ D^\mp_L)}\setminus N'$ the curve $\partial D'$ is homologous to $\partial D \pm
L',$ where $L'$ is a parallel copy of $L$ lying on $\partial N.$ Thus $M_{(\Sigma,\phi\circ
D^\mp_L)}$ is obtained from $M_{(\Sigma,\phi)}$ by a $\pm 1$ Dehn surgery on $L.$

We now must see that the contact structure one gets from $\pm 1$-contact surgery on $L$ is the
contact structure supported by the open book $(\Sigma,\phi\circ D^\mp_L).$ To do this we consider
the definition of compatible open book involving the Reeb vector field.
\bhw
Think about trying to show that $\xi$ can be isotoped arbitrarily close to the pages of the open
book. Intuitively this is not too hard to see (but as usual, making a rigorous proof out of this
intuition is not so easy).
\ehw
Since $(\Sigma, \phi)$ is compatible with $\xi$ there is a Reeb vector field $X$ for $\xi$ such that
$X$ is positively transverse to the pages and tangent to the binding. Notice our neighborhood
$N=N_\Sigma\times[\frac12-\epsilon,\frac 12+\epsilon]$ is such that $X$ is transverse to all the
$N_\Sigma\times\{t\}$ and the flow of $X$ preserves these pages. As usual we will now consider a
model situation. Consider $(\R^3, \ker (dz-y\, dx))/\sim,$ where $(x,y,z)\sim(x+1,y,z).$ It is easy
to arrange that the foliation on $N_\Sigma$ is the same as the foliation on $\{(x,y,z)| -\delta\leq
y\leq\delta\}.$ Thus we can assume the contact structure on $N$ is contactomorphic to the contact
structure on $N_m=\{(x,y,z)| -\epsilon\leq z\leq\epsilon, -\delta\leq y\leq\delta\}.$ Moreover, we
can assume this contactomorphism takes a Reeb vector field for $\xi$ to $\frac{\partial}{\partial
z},$ the Reeb vector field for $dz-y\, dx.$ (We do this by picking an identification of $N_\Sigma$
and the annulus in the $xy$-plane that preserves the characteristic foliation and then using the
Reeb vector fields to extend this identification.) Let $\psi:\partial N_m\to \partial N$ be the
diffeomorphism that agrees with the above contactomorphism everywhere except on $S_u=\{(x,y,z)|
z=\epsilon, -\delta\leq y\leq \delta\}\subset \partial N_m,$ where it differs by $D^\pm_L.$ Note that
gluing $N_m$ to $M\setminus N$ using $\psi$ will yield the manifold $M_{(L,\pm 1)}.$ Use $\psi$ to
pullback the characteristic foliation on $\partial(M\setminus N)$ to $\partial N_m.$ Note that the
characteristic foliation agrees with the characteristic foliation on $\partial N_m\setminus S_u$
induced by $\ker(dz-y\, dx).$
\bhw
Show there is a function $f:\{(x,y)| -\delta\leq y\leq \delta\}\to \R$ that equals $\epsilon$
near $|y|=\delta$ such that the characteristic foliation on the graph of $f$ agrees with the
foliation on $S_u$ that is pulled back from $N$ by $\psi.$ (By ``agrees with'' I mean that
$S_u$ and the graph of $f$ are isotopic rel boundary so that the isotopy takes the 
pulled back foliation to the characteristic foliation on the graph of $f.$) \hfill\break
HINT: Figure out what the pullback foliation is first. Then experiment with perturbing the
graph of the constant function.
\ehw
Now let $N_m'$ be the region bounded by $\partial N_m\setminus S_u$ union the graph of $f.$ There is
a natural way to think of $\psi$ as a map from $\partial N_m'$ to $\partial (M\setminus N)$ that
preserves the characteristic foliation. Furthermore we can extend $\psi$ to a neighborhood of the
boundary so that it preserves Reeb vector fields. The contact structures $\xi_{M\setminus N}$ and
$\ker(dz -y\, dx)|_{N_m'}$ glue to give a contact structure on $M_{(L,\pm 1)}.$ We can also glue up
the Reeb vector fields to get a Reeb vector field on $M_{(L,\pm 1)}$ that is transverse to the pages
of the obvious open book and tangent to the binding.
\end{proof}
An easy corollary of this theorem is the following.
\begin{thm}[Giroux 2002, \cite{Giroux??}]\label{steinfill}
A contact manifold $(M,\xi)$ is Stein fillable if and only if there is an open book decomposition
for $(M,\xi)$ whose monodromy can be written as a composition of right-handed Dehn twists.
\end{thm}
\begin{proof}
We start by assuming that there is an open book $(\Sigma, \phi)$ supporting $(M,\xi)$ for which
$\phi$ is a composition of right-handed Dehn twists. Let us begin by assuming that $\phi$ is the
identity map on $\Sigma.$ In Exercise~\ref{trivmongivesconsum} you verified that
$M=\#_{2g+n-1}S^1\times S^2,$ where $g$ is the genus of $\Sigma$ and $n$ is the number of boundary
components. Eliashberg has shown that $\#_{2g+n-1} S^1\times S^2$ has a unique strong symplectic
filling \cite{Eliashberg90a}. This filling $(W,\omega)$ is also a Stein filling. Thus we are done if
$\phi$ is the identity map.

Now assume $\phi=D^+_\gamma$ where $\gamma$ is a simple non-separating closed curve on $\Sigma.$ We
can use the Legendrian realization principle to make $\gamma$ a Legendrian arc on a page of the open
book. (Recall that even though our convex surface does not have Legendrian boundary we can still use the
Legendrian realization principle. See Remark~\ref{extralerp}.) (Note that we required $\gamma$ to be
non-separating so that we could use the Legendrian realization principle.) We know $(M_{(\Sigma,
id)},\xi_{(\Sigma,id)})$ is Stein filled by $(W,\omega)$ so by Theorems~\ref{surgeryandcobord} and
\ref{surgeryonbook} we can attach a 2-handle to $W$ to get a Stein filling of
$(M_{(\Sigma, D^+_\gamma)}, \xi_{(\Sigma, D^+_\gamma)}).$ If $\phi$ is a composition of more than
one right-handed Dehn twist along non-separating curves in $\Sigma$ we may clearly continue this
process to obtain a Stein filling of $(M_{(\Sigma, \phi)},\xi_{(\Sigma, \phi)}).$ The only thing
left to consider is when one or more of the curves on which we Dehn twist is separating. Suppose
$\gamma$ is separating. If both components of the complement of $\gamma$ contain parts of the
boundary of the page then we can still realize $\gamma.$ Thus we only have a problem when there is a
subsurface $\Sigma'$ of $\Sigma$ such that $\partial \Sigma'=\gamma.$ In this case we can use the
``chain relation'' (see Theorem~\ref{chainrel} in the Appendix) to write $D^+_\gamma$ as a
composition of positive Dehn twists along non-separating curves in $\Sigma'.$

For the other implication we assume that $(M,\xi)$ is Stein fillable by say $(W,J,\psi).$ According
to Eliashberg's Theorem~\ref{stein}, $W$ has a handle decomposition with only 1-handles and
2-handles attached along Legendrian knots with framing one less than the contact framing. Let $W'$
be the union of the 0- and 1-handles. Clearly $M'=\partial W'= \#_k S^1\times S^2$ and the induced
contact structure $\xi'$ is tight. So we have a Legendrian link $L$ in $M'$ on which we can perform
Legendrian surgery to obtain $(M,\xi).$ Now according to Corollary~\ref{Lonpage} there is an open
book decomposition $(\Sigma,\phi)$ for $(M',\xi')$ such that $L$ sits on a page of the open book. At
the moment $\phi$ might not be the composition of right-handed Dehn twists. According to a theorem
of Eliashberg \cite{Eliashberg92} and Colin \cite{Colin97} $M'$ has a unique tight contact
structure. We know there is a surface $\Sigma'$ such that $M'=M_{(\Sigma', id)}$ and the supported
contact structure is tight. According to Giroux's Stabilization Theorem~\ref{thm:stab} we can
positively stabilize $(\Sigma',id)$ and $(\Sigma,\phi)$ so that they become isotopic. Let
$(\Sigma'',\phi'')$ be their common stabilization. Since the stabilizations were all positive
$\phi''$ is a composition of positive Dehn twists. Moreover $L$ sits on a page of this open book. As
in the previous paragraph of the proof, performing Legendrian surgery on $L$ will change the open
book $(\Sigma'', \phi'')$ by composing $\phi''$ by right-handed Dehn twists. Thus we eventually get
an open book for $(M,\xi)$ whose monodromy consists of a composition of right-handed Dehn twists.
\end{proof}

\subsection{Proof of Theorem~\ref{thm:caps}}\label{subsec:capspf}
We are ready to begin the proof of Theorem~\ref{thm:caps}.
\begin{lem}[Etnyre 2004, \cite{Etnyre04}; {\em cf.} Stipsicz 2003, \cite{Stipsicz03}]\label{tostrong}
Suppose $(X,\omega)$ is a weak filling of $(M,\xi).$ Then there is a compact symplectic manifold
$(X',\omega')$ such that $(X,\omega)$ embeds in $(X', \omega')$ and the boundary of $(X',\omega')$
is strongly convex.
\end{lem}
\begin{proof}
Let $(B,\pi)$ be an open book for $(M,\xi).$ Using positive stabilizations of the open book we can
assume the binding is connected. Let $\phi$ be the monodromy of the open book. It is well known, see
Lemma~\ref{monoid} in the Appendix, that $\phi$ can be written
\[\phi=D_c^m\circ D_{\gamma_1}^{-1}\circ \ldots \circ D_{\gamma_n}^{-1},\]
where $\gamma_i$ are non-separating curves on the interior of the page $\Sigma$ and $c$ is a curve
on $\Sigma$ parallel to the boundary of $\Sigma.$

We know $M\setminus B$ is the mapping cylinder $\Sigma_\phi$ and an identical argument to the one in
the first paragraph of the proof of Theorem~\ref{surgeryonbook} says we may think of $\Sigma_\phi$
as
\[\coprod_{i=1}^{n} \Sigma_i/\sim,\]
 where $\Sigma_i=\Sigma\times[\frac{i-1}{n}, \frac{i}{n}]$ and $\sim$ is the equivalence relation
that glues $\Sigma\times\{\frac{i}{n}\}$ in $\Sigma_{i}$ to $\Sigma\times\{\frac{i}{n}\}$ in
$\Sigma_{i+1}$ by $D^{-1}_{\gamma_i}$ and $\Sigma\times{1}$ to $\Sigma\times{0}$ by $D_c^m.$ See
Figure~\ref{fig:monob}.
\begin{figure}[ht]
  \relabelbox \small {\epsfysize=1.3in\centerline{\epsfbox{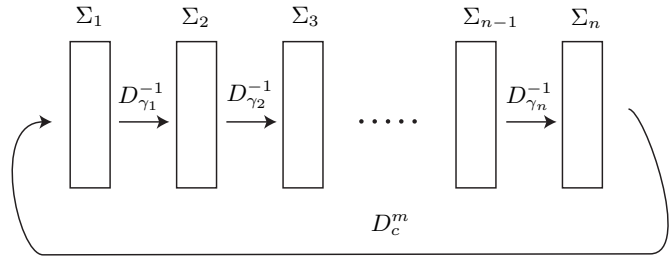}}}
  \adjustrelabel <-2pt,2pt> {1}{$D^{-1}_{\gamma_1}$} 
  \adjustrelabel <-2pt,2pt> {2}{$D^{-1}_{\gamma_2}$} 
  \adjustrelabel <-2pt,2pt> {n}{$D^{-1}_{\gamma_n}$}
  \relabel{m}{$D_c^m$} 
  \relabel{a}{$\Sigma_1$} 
  \relabel{b}{$\Sigma_2$} 
  \relabel{c}{$\Sigma_3$}
  \relabel{y}{$\Sigma_{n-1}$} 
  \relabel{z}{$\Sigma_n$} 
  \endrelabelbox
        \caption{Breaking up the monodromy.}
        \label{fig:monob}
\end{figure}

Choose a point $p\in(0, \frac{1}{n}).$ We can Legendrian realize $\gamma_1$ on the surface
$\Sigma\times\{p\}$ in $\Sigma_1\subset M.$ If we cut $\Sigma_\phi$ along $\Sigma\times\{p\}$ and
reglue using the diffeomorphism $D_{\gamma_1}$ then the new open book decomposition will have page
$\Sigma$ and monodromy $D_c^m\circ D_{\gamma_2}^{-1}\circ \ldots \circ D_{\gamma_n}^{-1}.$
Continuing in this way we can get an open book with page $\Sigma$ and monodromy $D_c^m.$ Denote the
contact manifold supported by this open book by $(M',\xi').$ By Theorem~\ref{surgeryonbook} we know
we can get from $(M,\xi)$ to $(M',\xi')$ by a sequence of Legendrian surgeries. Thus by
Theorem~\ref{surgeryandcobord} we can add 2-handles to $(X,\omega)$ in a symplectic way to get a
symplectic manifold $(X'',\omega'')$ with weakly convex boundary equal to $(M',\xi').$

Let $a_1,\ldots, a_{2g}$ be the curves on $\Sigma$ pictured in Figure~\ref{fig:basiscurves}.
\begin{figure}[ht]
  \relabelbox \small {\epsfysize=1.2in\centerline{\epsfbox{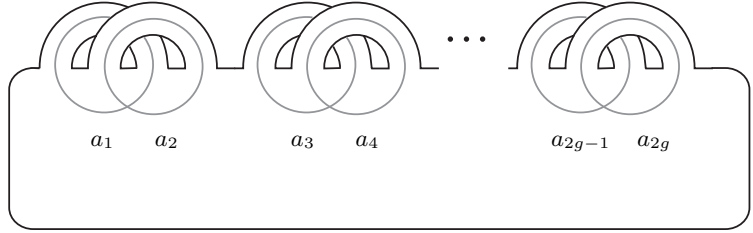}}} \relabel{1}{$a_1$}
  \relabel{2}{$a_2$} \relabel{3}{$a_3$} \relabel{4}{$a_4$} \relabel{5}{$a_{2g-1}$}
  \relabel{6}{$a_{2g}$} \endrelabelbox
        \caption{The curves $a_1,\ldots a_{2g}$ on $\Sigma.$}
        \label{fig:basiscurves}
\end{figure}
We can Legendrian realize these curves on separate pages of the open book for $(M',\xi')$ and do
Legendrian surgery on them to get the contact manifold $(M'',\xi'').$ Moreover, we can add 2-handles
to $(X'',\omega'')$ to obtain the symplectic manifold $(X',\omega')$ with weakly convex boundary
$(M'',\xi'').$ The open book supporting $(M'',\xi'')$ has page $\Sigma$ and monodromy $D_{a_1}\circ
\ldots \circ D_{a_{2g}}\circ D_c^m.$
\bhw
Show $M''$ is topologically obtained from $S^3$ by $\frac{1}{m}$ Dehn surgery on the knot in
Figure~\ref{fig:spic}.
\begin{figure}[ht]
  \relabelbox \small {\epsfysize=1in\centerline{\epsfbox{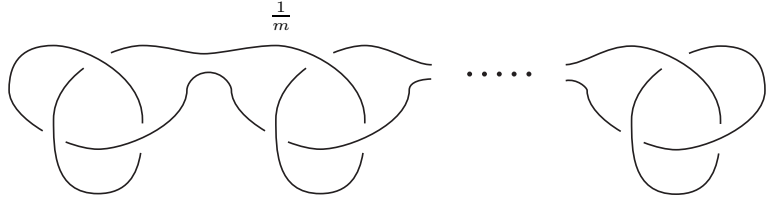}}}
  \relabel{m}{$\frac{1}{m}$} \endrelabelbox
        \caption{Topological description of $M''.$}
        \label{fig:spic}
\end{figure}
Thus $M''$ is a homology sphere.
\ehw
Now Theorem~\ref{filleqhom} says for a homology sphere a weak filling can be deformed into a strong
filling. Thus we may deform $(X',\omega')$ into a strong filling of $(M'',\xi'')$ and clearly
$(X,\omega)$ symplectically embeds into this manifold.
\end{proof}

\begin{lem}[Etnyre and Honda 2002, \cite{EtnyreHonda02a}; Gay 2002, \cite{Gay02}]\label{strongcap}
Given any contact manifold $(M,\xi)$ there is a strong concave filling of $(M,\xi).$
\end{lem}
This Lemma was also proven by Lisca and Mati\'c \cite{LiscaMatic97} for Stein fillable contact
structures. An alternate proof in the Stein case was provided by the work of Akbulut and Ozbagci
\cite{AO} coupled with that of Plamenevskaya \cite{Plam}. The proof below is in the spirit of Gay's
work.
\begin{proof}
We start with the symplectic manifold $(W,\omega_W)=(M\times[0,1], d(e^t\alpha)),$ where $\alpha$ is
a contact form of $\xi$ and $t$ is the coordinate on $[0,1].$ It is easy to see that $M\times\{0\}$
is a strongly concave boundary component of $W$ and $M\times\{1\}$ is a strongly convex boundary.
Our strategy will be to cap off the convex boundary component so we are left with only the concave
boundary component. Throughout this proof we will call the concave boundary of a symplectic manifold
the \dfn{lower boundary component} and the convex boundary component the \dfn{upper boundary
component}.

Let $(\Sigma,\phi)$ be an open book for $(M,\xi)$ with a connected boundary. As in the proof of
Lemma~\ref{tostrong} we can add 2-handles to $W$ to get a symplectic manifold $(W',\omega')$ with
lower boundary $M$ and upper boundary a contact manifold with open book having page $\Sigma$ and
monodromy $D^m_c$ where $c$ is a boundary parallel curve in $\Sigma.$

We want to argue that we can assume that $m=1.$ If $m<1$ then by adding symplectic 2-handles to
$(W',\omega')$ we can get to the situation where $m=1.$ (Note $c$ is separating, but we can handle
this as we did in the proof of Theorem~\ref{steinfill}.) Throughout the rest of the proof as we add
handles to $(W',\omega')$ we still denote the resulting manifold by $(W',\omega').$ We are left to
consider the situation where $m>1.$ For this we observe that we can increase the genus of $\Sigma$ as
follows. 
\bhw 
Show that if we add a symplectic 1-handle to $(W',\omega')$ this has the effect on the upper
boundary of $W'$ of connect summing with the standard (unique tight) contact structure on $S^1\times
S^2.$
\ehw 
But we know connect summing the contact manifold can be achieved by Murasugi summing their open
books. The tight contact structure on $S^1\times S^2$ has open book with page an annulus and
monodromy the identity map. Thus adding symplectic 1-handles to $W'$ has the effect on the open book
of the upper boundary component of adding a 1-handle to $\Sigma$ and extending the old monodromy
over this handle by the identity. So by adding 1-handles to $W$ we can arrange that the open book
for the upper boundary component has page $\Sigma'$ shown in Figure~\ref{fig:incgen}
\begin{figure}[ht]
  \relabelbox \small {\epsfysize=1.6in\centerline{\epsfbox{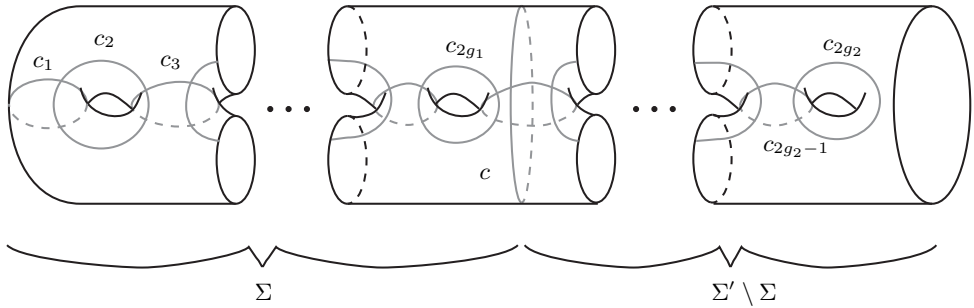}}} \relabel{1}{$c_1$}
  \relabel{2}{$c_2$} \relabel{3}{$c_3$} \relabel{4}{$c_{2g_1}$} \relabel{5}{$c_{2g_2-1}$}
  \relabel{6}{$c_{2g_2}$} \relabel{c}{$c$} \relabel{S}{$\Sigma$} \relabel{S2}{$\Sigma'\setminus
  \Sigma$} \endrelabelbox
        \caption{The surfaces $\Sigma$ and $\Sigma'$ and the curves $c$ and  $c_1,\ldots a_{2g_2}.$}
        \label{fig:incgen}
\end{figure}
and monodromy $D_c^m.$ Let $c_1,\ldots, c_{2g_1}$ be the curves in $\Sigma$ shown in
Figure~\ref{fig:incgen}; $g_1$ is the genus of $\Sigma.$ The Chain Relation (see
Theorem~\ref{chainrel}) says $D^m_c=(D_{c_1}\circ\ldots\circ D_{c_{2g_1}})^{m(4g_1+2)}.$ Now let
$c_{g_1+1},\ldots, c_{2g_2}$ be the curves shown in Figure~\ref{fig:incgen}; $g_2$ is the genus of
$\Sigma'.$
\bhw
Show that we can assume $m$ is such that we can choose the genus $g_2$ of $\Sigma'$ so that
$m(4g_1+2)=4g_2+2.$\hfill\break
HINT: Attach symplectic 2-handles.
\ehw
Thus we can attach symplectic 2-handles to $(W',\omega')$ so that the upper boundary $(M',\xi')$ has an
open book decomposition with page $\Sigma'$ and monodromy $D_{c'}=(D_{c_1}\circ\ldots\circ
D_{c_{2g_2}})^{4g_2+2},$ where $c'$ is a curve on $\Sigma'$ parallel to the boundary.
\bhw
Show $M'$ is an $S^1$ bundle over $\Sigma''$ with Euler number $-1,$ where $\Sigma''$ is $\Sigma'$ 
with a disk capping of its boundary.
\ehw
\bhw
Let $C$ be the $D^2$ bundle over $\Sigma''$ with Euler number 1. Construct a natural symplectic structure 
$\omega_C$ on $C.$ \hfill\break
HINT: On the circle bundle $\partial C$ there is a connection 1-form $\alpha$ that is also a contact
form on $\partial C.$ Use this to construct the symplectic form on $C.$ Note that if you think about the
symplectization of a contact structure you can easily get a symplectic structure on $C$ minus the
zero section. Some care is needed to extend over the zero section.
\ehw
\bhw
Show $(C,\omega_C)$ has a strongly concave boundary, $\partial C=-M',$ and $\xi'$ is the induced
contact structure.\hfill\break
HINT: The contact structure induced on $\partial C$ is transverse to the circle fibers. If you
remove a neighborhood of one of the fibers the resulting manifold is $\Sigma'\times S^1.$ This is
the mapping cylinder of the identity on $\Sigma'.$ Show that the contact planes can be isotoped
arbitrarily close to the pages. Now consider how the neighborhood of the fiber is glued back in.
\ehw
We now simply glue $(C,\omega_C)$ to the top of $(W',\omega')$ to get our concave filling of 
$(M,\xi).$
\end{proof}

We are now ready to prove the main result of this section.
\begin{proof}[Proof of Theorem~\ref{thm:caps}]
We start with a weak symplectic filling $(X,\omega)$ of $(M,\xi).$ Now apply Lemma~\ref{tostrong} to
embed $(X,\omega)$ symplectically into $(X',\omega')$ where $(X',\omega')$ has a strongly convex
boundary $(M',\xi').$ Now use Lemma~\ref{strongcap} to find a symplectic manifold $(X'',\omega'')$
that is a strong concave filling of $(M',\xi').$ Using an exercise from
Subsection~\ref{subsec:sfill} we can glue $(X',\omega')$ and $(X'',\omega'')$ together to get a
closed symplectic manifold $(W,\omega_W)$ into which $(X,\omega)$ embeds.
\end{proof}

\subsection{Sobering arcs and overtwisted contact structures}
Theorem~\ref{steinfill} gives a nice characterization of Stein fillable contact structures in terms
of open book decompositions. It turns out there is a similar characterization of overtwisted contact
structures due to Goodman \cite{Goodman}. Suppose we are given an oriented surface $\Sigma.$ Given
two properly embedded oriented arcs $a$ and $b$ on $\Sigma$ with $\partial a=\partial b$ we can
isotop them relative to the boundary so that the number of intersection points between the arcs is
minimized. At a boundary point $x$ of $a$ define $\epsilon (x)$ to be $+1$ if the oriented tangent
to $a$ at $x$ followed by the oriented tangent to $b$ at $x$ is an oriented basis for $\Sigma,$
otherwise we set $\epsilon (x)=-1.$ Let $i(a,b)=\frac12 (\epsilon(x) + \epsilon(y)),$ where $x$ and
$y$ are the boundary points of $a.$
\begin{defn}
Let $\Sigma$ be an oriented surface and $\phi:\Sigma\to \Sigma$ a diffeomorphism that fixes the
boundary. An arc $b$ properly embedded in $\Sigma$ is a \dfn{sobering arc} for the pair
$(\Sigma,\phi)$ if $i(b,\phi(b))\geq 0$ and there are no positive intersection points of $b$ with
$\phi(b)$ (after isotoping to minimize the number of intersection points).
\end{defn}
Note the definition of sobering arc does not depend on an orientation on $b.$
\begin{thm}[Goodman 2004, \cite{Goodman}]
If $(\Sigma,\phi)$ admits a sobering arc then the corresponding contact structure
$\xi_{(\Sigma,\phi)}$ is overtwisted.\qed
\end{thm}
We will not prove this theorem but indicate by an example how one shows a contact structure is
overtwisted if a supporting open book admits a sobering arc. Indeed, consider $(A,\phi)$ where
$A=S^1\times[-1,1]$ and $\phi$ is a left-handed Dehn twist about $S^1\times\{0\}.$ Of course this
is the open book describing the negative Hopf link in $S^3.$ Earlier we claimed that the associated
contact structure is overtwisted; we now find the overtwisted disk. (Actually we find a disk whose
Legendrian boundary violates the Bennequin inequality, but from this one can easily locate an
overtwisted disk.) The arc $b=\{pt\}\times[-1,1]\subset A$ is obviously a sobering arc. Let $T$ be
the union of two pages of the open book. Clearly $T$ is a torus that separates $S^3$ into two solid
tori $V_0$ and $V_1.$ We can think of $V_0$ as $A\times[0,\frac 12]$ (union part of the neighborhood
of the binding if you want to be precise) and $V_1$ is $A\times[\frac 12, 1].$ In $V_0$ we can take
$D_0=b\times[0,\frac 12]$ to be the meridional disk and in $V_1$ we take $D_1=b\times[\frac 12,1]$
to be the meridional disk. If we think of $T$ as the boundary of $V_1$ then $\partial D_1$ is simply
a meridional curve, that is, a $(1,0)$ curve. Note that $\partial D_0$ does not naturally sit on $T.$
We must identify $A\times\{\frac 12\}\subset V_0$ with $A\times\{\frac 12\}\subset V_1$ using the
identity map and $A\times\{0\}\subset V_0$ with $A\times\{1\}\subset V_1$ via $\phi.$ Thus $\partial
D_1$ in $T$ is a $(-1,1)$ curve. In particular these two meridional curves intersect once on $T.$ We can
Legendrian realize a $(0,1)$ curve on $T.$ Call the Legendrian curve $L.$ Note since
$(0,1)=(1,0)+(-1,1),$ $L$ bounds a disk $D$ in $S^3.$ With respect to the framing induced by $T$ the
twisting of $L$ is 0.
\bhw
Show that the framing induced on $L$ by $T$ is one larger than the framing induced by $D.$ Thus 
$tw(\xi, D)=1.$
\ehw
So we see that $\partial D$ violates the Bennequin inequality and thus $\xi$ is overtwisted. 
\bhw
Starting with $D$ find an overtwisted disk.
\ehw
In general, in the proof of the above theorem you will not always be able to find an explicit
overtwisted disk, but in a manner similar to what we did above you will always be able to construct
a Legendrian knot bounding a surface that violates the Bennequin inequality.

Lastly we have the following theorem.
\begin{thm}[Goodman 2004, \cite{Goodman}]
A contact structure is overtwisted if and only if it is there is an open book decomposition
supporting the contact structure that admits a sobering arc.
\end{thm}
The if part of this theorem is the content of the previous theorem. The only if part follows from:
\bhw
Show that any overtwisted contact structure is supported by an open book that has been negatively
stabilized. Show that this implies there is a sobering arc.\hfill\break
HINT: You will need to use Eliashberg's classification of overtwisted contact structures by their
homotopy class of plain field \cite{Eliashberg89}. If you are having trouble you might want to
consult \cite{Etnyre??} or, of course, the original paper \cite{Goodman}.
\ehw

\section{Appendix}\label{app}
We recall several important facts about diffeomorphisms of surfaces. First, given an embedded curve
$\gamma$ in an oriented surface $\Sigma$ let $N=\gamma\times[0,1]$ be a (oriented) neighborhood of
the curve. We then define the \dfn{right-handed Dehn twists along $\gamma$,} denoted $D_\gamma,$ to
be the diffeomorphism of $\Sigma$ that is the identity on $\Sigma\setminus N$ and on $N$ is given by
$(\theta,t) \mapsto (\theta+2\pi t,t),$ where $\theta$ is the coordinate on $\gamma=S^1$ and $t$ is
the coordinate on $[0,1]$ and we have chosen the product structure so that $\frac{\partial}{\partial
\theta}, \frac{\partial}{\partial t}$ is an oriented basis for $N\subset \Sigma.$ (Note that to make
$D_\gamma$ a diffeomorphism one needs to ``smooth'' it near $\partial N.$)
\begin{figure}[ht]
  \relabelbox \small {\epsfysize=1in\centerline{\epsfbox{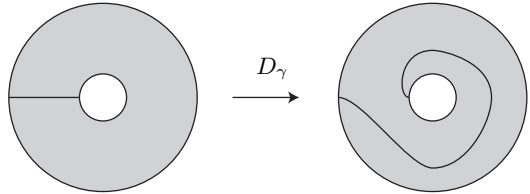}}}
  \relabel{D}{$D_\gamma$} \endrelabelbox
        \caption{A right-handed Dehn twist.}
        \label{fig:dt}
\end{figure}
A \dfn{left-handed Dehn twists about $\gamma$} is $D^{-1}_\gamma.$
\bhw
Show the following
\begin{enumerate}
\item $D_\gamma$ does not depend
on an orientation on $\gamma.$ 
\item If $\gamma$ and $\gamma'$ are isotopic then $D_\gamma$ and $D_{\gamma'}$ are isotopic 
  diffeomorphisms.
\end{enumerate}
\ehw

\begin{thm}[Lickorish 1962, \cite{Lickorish62}]
Any diffeomorphism of a compact oriented surface can be written as a composition of Dehn twists
about non-separating curves and curves parallel to the boundary of the surface.\hfill\qed
\end{thm}

There are several important relations among Dehn twists. For example:
\begin{enumerate}
\item For any $\gamma$ and diffeomorphism $f$ we have $f\circ D_\gamma\circ f^{-1}=D_{f(\gamma)}.$
\item If $\gamma$ and $\delta$ are disjoint then $D_\gamma\circ D_\delta=D_\delta\circ D_\gamma.$
\item If $\gamma$ and $\delta$ intersect in one point then $D_\delta\circ D_\gamma (\delta)$ is 
  isotopic to $\gamma.$
\item If $\gamma$ and $\delta$ intersect in one point then $D_\delta\circ D_\gamma\circ D_\delta = 
  D_\gamma\circ D_\delta \circ D_\gamma.$
\end{enumerate}
\bhw
Prove the above relations.
\ehw
\bhw
Show that given two non-separating curves $\gamma$ and $\delta$ on $\Sigma$ there is a 
diffeomorphism of $\Sigma$ taking $\gamma$ to $\delta.$
\ehw
In many of the above applications of open books we needed the following fundamental relation called 
the Chain Relation.
\begin{thm}\label{chainrel}
Let $\gamma_1,\ldots, \gamma_k$ be a chain of simple closed curves in $\Sigma,$ that is, the curves
satisfy $\gamma_i\cdot \gamma_j$ is $1$ if $|i-j|=1$ and is 0 otherwise, where $\cdot$ means
geometric intersection. Let $N$ be a neighborhood of the union of the $\gamma_i$'s. If $k$ is odd
then $N$ has two boundary components $d_1$ and $d_2.$ If $k$ is even then $N$ has one boundary
component $d.$ We have the following relations
\[
(D_{\gamma_1}\circ \ldots \circ D_{\gamma_k})^{2k+2} = D_d,\quad  \text{if $k$ is even, and}
\]
\[
(D_{\gamma_1}\circ \ldots \circ D_{\gamma_k})^{k+1} = D_{d_1}\circ D_{d_2}\quad \text{if $k$ is odd.}
\]\qed
\end{thm}
\bhw
Try to prove this theorem. Note: for $k=1$ it is trivial, for $k=2,3$ it is quite easy to explicitly 
check the relation.
\ehw
An important consequence of this theorem is the following lemma.
\begin{lem}\label{monoid}
Let $\Sigma$ be a surface with one boundary component, then any diffeomorphism of $\Sigma$ can be
written as the composition of right-handed Dehn twists about non-separating curves on the interior
of $\Sigma$ and arbitrary Dehn twists about a curve parallel to the boundary of $\Sigma.$
\end{lem}
\begin{proof}
We find a chain of curves $\gamma_1,\ldots, \gamma_{2g}$ in $\Sigma$ such that $\Sigma$ is a
neighborhood of their union. Thus the chain relation tells us that $(D_{\gamma_1}\circ \ldots \circ
D_{\gamma_{2g}})^{4g+2} = D_d.$ So clearly we can replace $D_{\gamma_i}^{-1}$ by a composition of
right-handed Dehn twists and one left-handed Dehn twist about $d.$ Now by the exercises above any
left-handed Dehn twist about a separating curve can be written as right-handed Dehn twists about
non-separating curves and a left-handed Dehn twist about $d.$
\end{proof}


\end{document}